\documentclass{elsart}

\usepackage{xspace}
\usepackage{amssymb}
\usepackage{array}
\usepackage{epsfig,psfrag}
\usepackage{version}

\makeatletter

\theoremstyle{plain}

\newtheorem{theorem}{Theorem}

\newtheorem{lemma}[theorem]{Lemma}

\newcommand {\be}{\begin{equation}}
\newcommand {\ee}{\end{equation}}
\newcommand {\bens}{\begin{eqnarray*}}
\newcommand {\eens}{\end{eqnarray*}}
\newcommand {\ben}{\begin{eqnarray}}
\newcommand {\een}{\end{eqnarray}}
\newcommand {\ba}{\begin{array}}
\newcommand {\ea}{\end{array}}
\newcommand {\bi}{\begin{itemize}}
\newcommand {\ei}{\end{itemize}}
\newenvironment{proof}[1][Proof]{\noindent\textbf{#1.} }{\ \rule{0.5em}{0.5em}}

\def\Rset{{\mathbb{R}}}
\def\Nset{{\mathbb{N}}}
\def\d{\mbox{d}}
\def\L{{\mathcal L}}
\def\kl{{\mathrm{kl}}}
\def\tv{{\mathrm{tv}}}
\def\egald{\mbox{$\,\stackrel{\mbox{\tiny d}}{=}\,$}}
\newcommand {\ie}{{\it i.e.}\xspace}
\newcommand {\eg}{{\it e.g.}\xspace}

\begin{document}

\begin{frontmatter}

  \title{On classes of  non-Gaussian asymptotic minimizers in entropic
    uncertainty principles}

  \author[LIS]{S.          Zozor}        and        \author[LIS,UM]{C.
    Vignat}

  \address[LIS]{Laboratoire  des  Images  et  des Signaux\\Rue  de  la
    Houille
    Blanche\\B.P. 46\\38420 Saint Martin d'H\`eres Cedex\\France\\
    Phone \# +33 4 76 82 64 23\\
    Fax \# +33 4 76 82 63 84\\
    E-mail: {\tt steeve.zozor@lis.inpg.fr}}

  \address[UM]{Institut   Gaspard  Monge\\Universit\'e  de   Marne  la
    Vall\'ee\\77454  Marne  La  Vall\'ee  Cedex\\France\\E-mail:  {\tt
      christophe.vignat@lis.inpg.fr}}

\begin{abstract}
  In  this  paper  we   revisit  the  Bialynicki-Birula  \&  Mycielski
  uncertainty  principle \cite{BiaMyc75}  and its  cases  of equality.
  This   Shannon  entropic  version   of  the   well-known  Heisenberg
  uncertainty principle  can be used when dealing  with variables that
  admit  no  variance.  In  this  paper,  we  extend this  uncertainty
  principle to R\'enyi  entropies. We recall that in  both Shannon and
  R\'enyi  cases, and  for a  given dimension  $n$, the  only  case of
  equality occurs  for Gaussian random  vectors.  We show that  as $n$
  grows,  however, the bound  is also  asymptotically attained  in the
  cases of  $n$-dimensional Student-t and  Student-r distributions.  A
  complete  analytical study  is  performed  in a  special  case of  a
  Student-t distribution.   We also show numerically  that this effect
  exists for the particular case of a $n$-dimensional Cauchy variable,
  whatever the  R\'enyi entropy  considered, extending the  results of
  Abe \cite{AbeRaj01} and illustrating the analytical asymptotic study
  of the student-t  case.  In the Student-r case,  we show numerically
  that  the same  behavior occurs  for uniformly  distributed vectors.
  These particular cases and other ones investigated in this paper are
  interesting since they show  that this asymptotic behavior cannot be
  considered as a ``Gaussianization'' of the vector when the dimension
  increases.
\end{abstract}

\begin{keyword}
  Entropic uncertainty relation, R\'enyi/Shannon entropy, multivariate
  Student-t and Student-r distributions.

  \PACS 02.50.-r, 02.50.Sk, 65.40.Gr
\end{keyword}
\end{frontmatter}

%%%%%%%%%%%%%%%%%%%%%%%%%%%%%%%%%%%%%%%%%%%%%%%%%%%%%%%%%%%%%%%%%%%%%

\section{Introduction}
\label{introduction:sec}

Let us consider an  $n$-dimensional wave packet $\Psi_n(x)$ and denote
by   $\widehat{\Psi}_n(u)   =   (2   \pi)^{-\frac{n}{2}}\int_{\Rset^n}
\Psi_n(x) \mbox{e}^{-  \imath u^t  x} \, \d  x$ its  Fourier transform
(consider for example the position of a particle and its momentum). In
the following, we will denote  by $X_n$ a zero-mean random vector with
probability  density function  (pdf) $f_n(x)  = |\Psi_n(x)|^2$  and by
$\widetilde{X}_n$  a  random  vector  with pdf  $\widetilde{f}_n(x)  =
|\widehat{\Psi}(x)|^2$  (by  Parseval's  relation,  this  is  a  pdf).
Vectors  $X_n$  and $\widetilde{X}_n$  are  called conjugated.   The
well-known  Heisenberg  uncertainty  principle (H.U.P.)   relates  the
``information''  available in two  conjugated random  vectors, stating
that  the product of  their variances  is larger  than a  given bound,
namely

\be
\frac{\left( E \left[ X_n^t  X_n \right] \, E \left[ \widetilde{X}_n^t
      \widetilde{X}_n     \right]     \right)^{\frac{1}{2}}}{n}    \ge
\frac{1}{2}
\label{Heisenberg:eq}
\ee
(see  also  \cite{DemCov91} for  a  matrix-variate  extension of  this
result).  The H.U.P.   is important in physics since  it expresses the
impossibility of an arbitrarily  accurate preparation of both position
and momentum of a particle.  This inequality finds also application in
areas   of   signal  processing   such   as  time-frequency   analysis
\cite{Fla99,Fla01}   where  it  is   known  as   the  Heisenberg-Gabor
inequality.  However, the H.U.P. has  a meaning only if the quantities
in balance exist.

To bypass  this restriction, Bialynicki-Birula \&  Mycielski showed in
1975 \cite{BiaMyc75}  that the H.U.P.  can be  extended to information
theoretic measures:  more precisely, they  showed that the sum  of the
Shannon entropy  rates of $X_n$ and of  $\widetilde{X}_n$ verifies the
Bialynicki-Birula \& Mycielski inequality (B.B.M.I.)
\be
\frac{H(X_n) + H(\widetilde{X}_n)}{n} \ge 1 + \log \pi
\label{uncertainty_Shannon:eq}
\ee
where the  Shannon entropy  is $H(X_n) =  - \displaystyle \int  f_n \,
\log f_n$ (and likewise  for $\widetilde{X}_n$).  The B.B.M.I.  can be
expressed equivalently via the  entropy power $N(X_n) = \frac{1}{2 \pi
  e}  \,  \exp  \left(   \frac{2}{n}  H(X_n)  \right)$  (as  given  in
\cite{DemCov91}) by
\be
\left(    N(X_n)    N(\widetilde{X}_n)    \right)^{\frac{1}{2}}    \ge
\frac{1}{2}.
\ee

The  B.B.M.I.   (\ref{uncertainty_Shannon:eq})  is stronger  than  the
H.U.P. (\ref{Heisenberg:eq}):  as shown in \cite{BiaMyc75},  it can be
applied also to variables  with infinite variance, provided that their
Shannon entropy  exists: it  is the  case for a  Cauchy pdf  $f_n(x) =
\frac{\Gamma \left( \frac{n+1}{2}  \right)}{\pi^{\frac{n+1}{2}} } (1 +
x^t  x)^{-\frac{n+1}{2}}$,  for example.   Moreover,  it  is shown  in
\cite{BiaMyc75} that  B.B.M.I.  (\ref{uncertainty_Shannon:eq}) implies
H.U.P.  (\ref{Heisenberg:eq}) when dealing with variables that admit a
variance\footnote{Since under covariance  constraint $N$ is maximum in
  the Gaussian context,  $N(X_n) \le N(G_n) = n  \sigma^2$ where $G_n$
  is Gaussian with  the same covariance than $X_n$,  and similarly for
  $\widetilde{X}_n$: the  product of the variances is  higher than the
  product of  the entropy power, implying the  H.U.P.}.  As inequality
(\ref{Heisenberg:eq}),  the  B.B.M.I.   finds applications  in  signal
processing  (see  \cite{Fla99,Fla01}   and  references  therein).   In
physics, maximization  without constraint of the sum  of the entropies
that appear in the B.B.M.I.   has been suggested \cite{GadBen85} as an
interesting  counterpart  to  the  classical ``maximum  entropy  under
constraint''  approach  for  the  derivation  of  the  wave  functions
associated with atomic systems.  As for the Heisenberg inequality, the
lower  bound in  (\ref{uncertainty_Shannon:eq}) is  attained  only for
Gaussian wave  packets.  But in their paper  \cite{AbeRaj01}, Abe {\it
  et al.}   showed that this  bound is also asymptotically  reached by
$n-$dimensional Cauchy vectors when  $n$ increases.  In this paper, we
will  extend   this  observation  by  exhibiting   other  families  of
distributions that  show the same  asymptotical behavior.  Moreover,we
will focus on the possible  interpretations of this behavior, namely a
``Gaussianization effect'':  one natural interpretation  is that these
distributions get closer to the  Gaussian distribution - in some sense
to  be determined -  as the  dimension increases;  this interpretation
will  be  proved  erroneous  by   providing  some  cases  in  which  a
distribution  reaches asymptotically  the bound  of the  B.B.M.I.  but
keeps   at  non-zero   (even  infinite)   distance  of   any  Gaussian
distribution.   As  we will  explain  in  the  conclusion, the  effect
observed is  mainly due to the  normalization $1/n$, which  may be too
strong in the considered non-iid context.

This paper is organized as follows:
\begin{itemize}
\item  In the first  part, we  come back  to the  Bialynicki-Birula \&
  Mycielski uncertainty relation, we  reformulate it and we generalize
  it to  the R\'enyi  entropies.  We then  give the expression  of the
  entropies  rates under  consideration when  dealing  with elliptical
  random vectors.
\item We then address the study of the asymptotic cases of equality in
  the   generalized  B.B.M.I.,  particularizing   to  the   family  of
  $n$-dimensional  Student-t variables  with $m$  degrees  of freedom.
  For this class of variables, we  provide an upperbound of the sum of
  the  entropy  rates  in   balance  which  permits  to  evaluate  the
  asymptotic behavior of this  quantity.  We then perform the complete
  analytical study in the particular case $m=n+2$, confirming that the
  lower bound  of the uncertainty relation  is attained asymptotically
  as $n$ increases.   In a second illustration, we  revisit the Cauchy
  case  ($m=1$)  as studied  by  Abe in  the  context  of the  R\'enyi
  formulation of the uncertainty relation, and we show that the effect
  observed by Abe remains for any ``admissible'' R\'enyi entropy.

  Secondly, we  explore what happens in  the $n$-dimensional Student-r
  class of variables with $m$ degrees of freedom.  We present the case
  $m=n$,   corresponding   to   the   uniform  distribution   in   the
  $n$-dimensional sphere, as well as some other cases.
\item Finally,  we provide  some clarifications about  this asymptotic
  behavior, explaining  why a Gaussianization effect,  measured in the
  distribution  or  in  the   information  divergence  sense,  is  not
  necessary to reach asymptotically equality in the B.B.M.I.
\end{itemize}

%%%%%%%%%%%%%%%%%%%%%%%%%%%%%%%%%%%%%%%%%%%%%%%%%%%%%%%%%%%%%%%%%%%%%

\section{The R\'enyi entropy uncertainty relation}
\label{Renyi:sec}

The proof  of the B.B.M.I. (\ref{uncertainty_Shannon:eq})  is based on
the  Beckner inequality  relating  the norms  of  any (wave)  function
$\Psi_n$   of   $L^p(\Rset^n)$    and   of   its   Fourier   transform
$\widehat{\Psi}_n$, \ie
\be
\|  \widehat{\Psi}_n \|_q \le (C_{p,q})^n \,  \| \Psi_n \|_p
\label{Beckner:eq}
\ee
where $p$ and $q$ are conjugated, \ie $\frac{1}{p} + \frac{1}{q} = 1$,
where $p  \in ]1  ; 2]$  and where $C_{p,q}$  is the  Babenko constant
expressed as  $C_{p,q} = \left( \frac{2  \pi}{p} \right)^{- \frac{1}{2
    p}}    \left(    \frac{2    \pi}{q}    \right)^{\frac{1}{2    q}}$
\cite{Bec75:02,Bec75}.  In  \cite{BiaMyc75}, the authors  consider the
positive   function  $W(q)  =   (C_{p,q})^n  \|   \Psi_n  \|_p   -  \|
\widehat{\Psi}_n  \|_q$ defined for  $q \ge  2$. Since  $W(2) =  0$ by
Parseval's identity  and since $W(q)$  is positive, the  derivative of
$W(q)$  in  $q=2$ is  positive  as  well,  which leads  to  inequality
(\ref{uncertainty_Shannon:eq}).

The B.B.M.I.  (\ref{uncertainty_Shannon:eq}) can be  extended to other
measures of information such as the R\'enyi entropies that include the
Shannon entropy as a special case.  The R\'enyi entropy with parameter
$\lambda$ is defined as
\be
H_\lambda(X_n) =  \frac{1}{1-\lambda} \, \log  \left( \int f_n^\lambda
\right)
\ee
for  $\lambda \ne  1$ \cite{Ren61}.   When  $\lambda$ tends  to 1,  by
l'Hospital's rule,  $H_\lambda$ converges to the  Shannon entropy that
will thus be denoted by continuity  $H_1 = H$.  The R\'enyi entropy is
widely used,  not only in physics (\eg  statistical mechanics, physics
of  turbulence, cosmology,  see  \cite{LenMen00,JizAri04,ParBir05} and
references  therein), but  in various  other areas  such as  in signal
processing (time scale  analysis, decision problems, machine learning,
see  \cite{Fla99,CovTho91,BarFla01,KreKas05} and  references therein),
or   image  processing   (image  matching,   image   registration  see
\cite{HerMa02,NeeHer05} and references therein).

\begin{theorem}
  With  $\frac{1}{p}+\frac{1}{q}=1$   and  for   any  $p  >   1$,  the
  B.B.M.I. writes in terms of R\'enyi entropy as
\be
\frac{H_{\frac{p}{2}}(X_n)  + H_{\frac{q}{2}}(\widetilde{X}_n)}{n} \ge
\log (2 \pi) + \frac{\log p}{p-2} + \frac{\log q}{q-2}.
\label{uncertainty_Renyi:eq}
\ee
\end{theorem}

\begin{proof}
It is straightforward that
\bens
\log  \|\Psi_n\|_p   =  \log  \left\|  f_n^{1/2}  \right\|_p   &  =  &
\frac{1}{p} \log \left( \int f_n^{p/2} \right)
\\
& = & \frac{2 - p}{2  p} \, H_{\frac{p}{2}}(X_n).
\eens
Hence,  taking  the logarithm  of  both  sides of  (\ref{Beckner:eq}), 
using the  conjugation relation $\frac{1}{p} +  \frac{1}{q} = 1$
and $1 < p \le 2$ leads to
\bens
\frac{2-p}{2  \, p} \,  H_{\frac{p}{2}}(X_n) -  \frac{2-q}{2 \,  q} \,
H_{\frac{q}{2}}(\widetilde{X}_n) &  \ge & n \left( \frac{1}{2  \, p} -
  \frac{1}{2 \, q} \right) \, \log(2 \pi)
\\
& & + \frac{n}{2 \, p} \log p - \frac{n}{2 \, q} \log q.
\eens
But  $p$  and  $q$  are  conjugated  so that  $\frac{2-q}{2  \,  q}  =
\frac{p-2}{2 \, p}$.  Hence, since  $2-p \ge 0$, multiplying each side
by  $\frac{2  \,  p}{n  \,  (2-p)}$ one  finally  obtains  uncertainty
relation    (\ref{uncertainty_Renyi:eq}).    Since    the    pdfs   of
$\widetilde{\widetilde{X}}_n$  and   of  $X_n$  coincide,   $X_n$  and
$\widetilde{X}_n$      have       a      symmetrical      role      in
(\ref{uncertainty_Renyi:eq})   and  can  then   be  exchanged:   as  a
consequence inequality (\ref{uncertainty_Renyi:eq}) holds for any $p >
1$, provided that the entropies in balance exist.
\end{proof}

By  taking the limit  $p \to  2$ in  (\ref{uncertainty_Renyi:eq}), the
B.B.M.I.  (\ref{uncertainty_Shannon:eq}) is recovered, proving that it
is a particular case  of (\ref{uncertainty_Renyi:eq}).  We notice that
a   similar   generalization    exists   for   the   Tsallis   entropy
\cite{MaaUff88,Raj95,GhoCha00},  which is  widely used  in statistical
physics  and   related  to  the  R\'enyi  entropy   by  an  invertible
transformation  \cite{LenMen00,JizAri04,ParBir05,TsaMen98,Tsa99}.   In
this paper  we will focus on  the R\'enyi entropy  since the resulting
form  of the  uncertainty relation  is  very similar  to the  B.B.M.I.
Furthermore, the quantity $\frac{1}{n} H_\lambda (X_n)$, also known as
the entropy rate (\ie entropy per sample), is very often considered in
the information theoretic context \cite{CovTho91}.

\begin{theorem}
  For a given $n$,  case of equality in (\ref{uncertainty_Shannon:eq})
  or (\ref{uncertainty_Renyi:eq}) is reached  if and only if $X_n$ and
  $\widetilde{X}_n$ are Gaussian random vectors.
\end{theorem}

\begin{proof}
  It  is straightforward  to  check that  Gaussian  vectors $X_n$  and
  $\widetilde{X}_n$    reach    equality    in    either    inequality
  (\ref{uncertainty_Shannon:eq})  or  (\ref{uncertainty_Renyi:eq})  by
  plugging the Gaussian  pdfs.  Conversely, for $p \ne  2$, since only
  Gaussian waves  functions achieve equality  in (\ref{Beckner:eq}) as
  it is proved in \cite{Lie90},  the Gaussian waves are the {\em only}
  wave packets which achieve equality in (\ref{uncertainty_Renyi:eq}).
  The  Shannon case $p=2$  is more  subtle, but  it has  been recently
  proved     that     equality     is    reached     in     inequality
  (\ref{uncertainty_Shannon:eq})  {\em  only}  in  the  Gaussian  case
  \cite{OzaPrz04}.
\end{proof}

However,  in  \cite{AbeRaj01},  Abe  showed  in the  case  of  Shannon
entropies   that  $n-$variate   Cauchy  vectors   reach   equality  in
(\ref{uncertainty_Shannon:eq}) asymptotically with  $n$.  We will show
in  the next  part  that this  result  extends in  fact to  inequality
(\ref{uncertainty_Renyi:eq}) for any value  $p/2$ of the entropy index
for  which the  entropy exists.   Furthermore  we will  show that  the
Cauchy case is  not the only one exhibiting  this behavior.

%%%%%%%%%%%%%%%%%%%%%%%%%%%%%%%%%%%%%%%%%%%%%%%%%%%%%%%%%%%%%%%%%%%%%

\section{Asymptotic cases of  equality}
\label{asymptotic_equality:sec}

As   previously    stated,   although   equality    is   achieved   in
(\ref{uncertainty_Renyi:eq})   only    by   Gaussian   wave   packets,
$n-$dimensional Cauchy wave packets reach equality asymptotically with
the dimension  $n$.  An important  question is to understand  what are
the ingredients that lead to this asymptotic behavior.

We first  remark that in the Cauchy  case presented by Abe  and in the
cases  studied  below,  the   components  of  the  random  vector  are
dependent: in other words the wave function $\Psi_n$ is not separable.
This  condition is  clearly required  since, dealing  with independent
identically distributed  (i.i.d.) components (separable  wave function
$\Psi_n$), the  entropy rate  coincides with the  entropy of  a single
component:  $H_\lambda(X_n)  =  n   H_\lambda(X)$  where  $X$  is  any
component of  $X_n$.  But  if $X_n$ is  i.i.d., the  conjugated vector
$\widetilde{X}_n$ is i.i.d.  as well  since the Fourier transform of a
separable  function  is   separable.   Hence  with  obvious  notations
$H_\lambda(\widetilde{X}_n)   =    n   H_\lambda(\widetilde{X})$   and
$\frac{H_\lambda(X_n) +  H_\lambda(\widetilde{X}_n)}{n} = H_\lambda(X)
+  H_\lambda(\widetilde{X})$: as a  consequence, no  asymptotic effect
can appear in the i.i.d. setup.

Furthermore,   for   any    invertible   matrix   $M$,   the   R\'enyi
$\lambda$-entropy of $Y_n  = M X_n$ is expressed  as $H_\lambda(Y_n) =
\log  |M|  +  H_\lambda(X_n)$  and  since  $\widetilde{Y}_n  =  M^{-t}
\widetilde{X}_n$,  the sum  of  the entropy  rates  is (matrix)  scale
invariant: it is thus impossible to reach equality by introducing such
simple  correlation  between the  components  of  the random  vectors.

Except for these basic requirements, the answer remains open as far as
we know.   The study  of the following  cases, showing  very different
behaviors, attempts to give some elements of answer.

%%%%%%%%%%%%%%%%%%%%%%%%%%%%%%%%%%%%%%%%%%%%%%%%%%%%%%%%%%%%%%%%%%%%%

\subsection{Derivation of the entropy rates in the elliptic case}
\label{elliptic:sec}

We concentrate in the following on the case of elliptical
random vectors.  A vector $X_n$  is called elliptical if its pdf $f_n$
is a single-valued function of a quadratic form \cite{Chu73}
\be
f_n(x) =  |\Sigma_n^{-1}|^{- \, \frac{1}{2}} \, d_n  \left( \left( x^t
    \Sigma_n^{-1} x \right)^{\frac{1}{2}} \right)
\label{fndn:eq}
\ee
for some  function $d_n$ and  where $\Sigma_n$ is a  positive definite
symmetric    matrix.
$\Sigma_n$ is called characteristic matrix.
In other words, the  random vector $\Sigma_n^{-1/2} X_n$ is isotropic.
Due  to   the  matrix  scale   invariance  of  the   studied  entropic
inequalities evoked above, we will consider without loss of generality
in  the following  that  $\Sigma_n$ is  proportional  to the  identity
matrix $I_n$.  Otherwise, except in (\ref{uncertainty_Renyi:eq}) where
there is no influence,  $X_n$ and $\widetilde{X}_n$ must be understood
as   $\Sigma_n^{-1/2}   X_n$   and  $\Sigma_n^{1/2}   \widetilde{X}_n$
respectively.

\begin{theorem}
  If $X_n$ is elliptical  as in (\ref{fndn:eq}) (with $\Sigma_n=I_n$),
  then the sum of the entropy rates of $X_n$ and $\widetilde{X}_n$ is
\be\ba{lll}
U_p(X_n)       &      =      &       \frac{H_{\frac{p}{2}}(X_n)      +
  H_{\frac{q}{2}}(\widetilde{X}_n)}{n}\vspace{2.5mm}
\\
& =  & \displaystyle \frac{2}{n}  \log \left(\frac{2 \pi^{n/2}}{\Gamma
    \left( n/2 \right)} \right)\vspace{2.5mm}
\\
&& \displaystyle + \frac{2}{n  (2-p)} \log \!\!  \int_0^{+\infty} \!\!
r^{\frac{(n-1)(2-p)}{2}} D_n(r)^{\frac{p}{2}} \d r\vspace{2.5mm}
\\
&& \displaystyle + \frac{2}{n  (2-q)} \log \!\!  \int_0^{+\infty} \!\!
r^{\frac{(n-1)(2-q)}{2}} E_n(r)^{\frac{q}{2}} \d r
\label{U_Renyi_si:eq}
\ea\ee
for $p \ne 2$ and
\ben
\lefteqn{U_2(X_n) = \frac{H(X_n) + H(\widetilde{X}_n)}{n}} \nonumber
\\
&&  =  \frac{2}{n}  \log  \left(\frac{2 \pi^{n/2}}{\Gamma  \left(  n/2
    \right)}  \right) +  \frac{n-1}{n} \,  \int_0^{+\infty} \log  r \,
(D_n(r) +
E_n(r)) \, \d r\nonumber
\\
&&  - \frac{1}{n}  \, \int_0^{+\infty}  (D_n(r) \log(D_n(r))  + E_n(r)
\log(E_n(r)) ) \, \d r
\label{U_Shannon_si:eq}
\een
for $p=2$, where
\be\left\{\ba{lll}
  D_n(r) & =  & \displaystyle \frac{2 \pi^{\frac{n}{2}}}{\Gamma \left(
      \frac{n}{2} \right)} \, r^{n-1} d_n(r)\vspace{5mm}
  \\
  E_n(r)   &   =  &   \displaystyle   \left(  \int_0^{+\infty}   (\rho
    r)^{\frac{1}{2}} D_n(\rho)^{\frac{1}{2}} J_{\frac{n}{2}-1}(\rho r)
    \, \d \rho \right)^2
\ea\right.
\label{DnEn:eq}
\ee
are the pdfs of  the Euclidean norms $\|X_n\|$ and $\|\widetilde{X}_n\|$
respectively.
\end{theorem}

\begin{proof}
  Using spherical coordinates, a simple computation shows that the pdf
  $D_n(r)$ of the Euclidean norm $\| X_n\|$ of $X_n$ is
\be
D_n(r) = \frac{2  \pi^{\frac{n}{2}}}{\Gamma \left( \frac{n}{2} \right)}
  r^{n-1} d_n(r)
\label{Dn:eq}
\ee
if $r \ge 0$ and is zero otherwise (see also \cite[eq. (7)]{Lor54}).

Let  us show now  that, $X_n$  being elliptical,  $\widetilde{X}_n$ is
elliptical:  the pdf  $\widetilde{f}_n$ of  the conjugate  variable is
given by
\bens
\widetilde{f}_n(u) & = & \left| (2 \pi)^{- \frac{n}{2}} \int_{\Rset^n}
  f_n^{\frac{1}{2}}(x) \mbox{e}^{-\imath u^t x} \d x \right|^2
\\
& = & \left|  (2 \pi)^{- \frac{n}{2}} \int_{\Rset^n} d_n^{\frac{1}{2}}
  \left((x^t  x)^{\frac{1}{2}} \right) \mbox{e}^{-\imath  u^t x}  \d x
\right|^2
\eens
Applying \cite[eq. (5)]{Lor54}, we obtain
\bens
\widetilde{f}_n(u)  = \left| \int_0^{+\infty}  \rho^{\frac{n}{2}} (u^t
  u)^{-  \frac{n-2}{4}} J_{\frac{n}{2}-1}(\rho  (u^t u)^{\frac{1}{2}})
  d_n^{\frac{1}{2}}(\rho) \d \rho \right|^2
\eens
where  $J_\nu$ is  the Bessel  function of  the first  kind  and order
$\nu$.  This result proves that $\widetilde{X}_n$ is elliptical.  From
the   above   expression    of   $\widetilde{f}_n$   and   the   forms
(\ref{fndn:eq})-(\ref{Dn:eq}), the  pdf $E_n$ of
$\|\widetilde{X}_n\|$ writes then
\be
E_n(r)    =    \left(    \int_0^{+\infty}    (\rho    r)^{\frac{1}{2}}
  D_n^{\frac{1}{2}}(\rho)   J_{\frac{n}{2}-1}(\rho  r)   \,   \d  \rho
\right)^2
\label{En:eq}
\ee
for $r \ge 0$ and zero  otherwise.  We remark that $E_n$ is the square
of the $(\frac{n}{2}-1)$-order Hankel transform of $D_n^{\frac{1}{2}}$
(see also \cite{Bat54,Pou00}).

Now the $\lambda$-R\'enyi entropy of $X_n$ is 
\bens
H_\lambda(X_n)  & =  & \frac{1}{1-\lambda}  \log \int_{\Rset^n}
f_n^\lambda(x) \, \d x
\\
&     =     &      \frac{1}{1-\lambda}     \log     \left(     \frac{2
    \pi^{\frac{n}{2}}}{\Gamma      \left(     \frac{n}{2}     \right)}
  \int_0^{+\infty} r^{n-1} d_n^\lambda(r) \, \d r \right)
\eens
the  last  equality being  an  application of  \cite[4.642]{GraRyz80}.
This yields the following expression for the $\lambda$-R\'enyi entropy
of $X_n$,
\ben
H_\lambda(X_n)    =    \log    \frac{2\pi^{n/2}}{\Gamma(\frac{n}{2})}+
\frac{1}{1-\lambda}    \log    \int_0^{+\infty}   r^{(n-1)(1-\lambda)}
D_n^\lambda(r) \, \d r
\label{Renyi_Xn:eq}
\een
The Shannon entropy of $X_n$ can be deduced from (\ref{Renyi_Xn:eq})
using L'Hospital's rule:
\ben
\hspace{-7mm} 
H(X_n)   =   \log   \frac{2\pi^{n/2}}{\Gamma(\frac{n}{2})}   +   (n-1)
\int_0^{+\infty} D_n(r) \log r \,  \d r - \int_0^{+\infty} D_n(r) \log
D_n(r) \, \d r
\label{Shannon_Xn:eq}
\een
We have shown  that $\widetilde{X}_n$ is again elliptical  so that its
$\lambda$-R\'enyi  entropy is  again expressed  as (\ref{Renyi_Xn:eq})
with $D_n$, the pdf of $\| X_n  \|$, replaced by $E_n$, the pdf of $\|
\widetilde{X}_n  \|$,  leading  to results  (\ref{U_Renyi_si:eq})  and
(\ref{U_Shannon_si:eq}).
\end{proof}

We have now  all the material to investigate  more deeply some special
cases of asymptotic equalities in (\ref{uncertainty_Renyi:eq}).

%%%%%%%%%%%%%%%%%%%%%%%%%%%%%%%%%%%%%%%%%%%%%%%%%%%%%%%%%%%%%%%%%%%%%

\subsection{The general Student-t case}
\label{Student-t:sec}

We consider  in this  subsection the case  where $X_n$  is distributed
according   to  a  Student-t   law  with   $m$  degrees   of  freedom,
\ie\footnote{More  rigorously the  random variable  $\sqrt{m}  X_n$ is
  Student-t  \cite{AbrSte70}, but  we have  seen that  the sum  of the
  entropy rates is insensitive to any scaling factor.}
\be
f_n(x) =  \frac{\Gamma \left( \frac{n+m}{2} \right)}{\pi^{\frac{n}{2}}
  \Gamma \left(  \frac{m}{2} \right)} \,  \left( 1 + x^t  x \right)^{-
  \frac{n+m}{2}}
\label{pdf_Student-t:eq}
\ee
where  $ m >  0$.  When  $m=1$, $X_n$  is the  well-known multivariate
Cauchy vector,  as studied by  Abe in \cite{AbeRaj01}.   The Student-t
variables  play  an important  role  in  statistics  because of  their
power-law  behavior and their  simple analytic  expression.  Moreover,
they maximize the R\'enyi/Tsallis  entropy with parameter $\lambda = 1
-        \frac{2}{n+m}$       under        covariance       constraint
\cite{ParBir05,TsaMen98,Tsa99,VigHer04,Bas04}.

\begin{theorem}
  If $X_n$ is Student-t distributed, then the pdf of $\| X_n \|$ is
\be
D_n(r)  = \frac{2 \Gamma  \left( \frac{n+m}{2}  \right)}{\Gamma \left(
    \frac{n}{2}   \right)  \Gamma   \left(  \frac{m}{2}   \right)}  \:
\frac{r^{n-1}}{(1+r^2)^{\frac{m+n}{2}}}
\label{Dn_Student-t:eq}
\ee
while the pdf of $\| \widetilde{X}_n \|$   is
\be
E_n(r)   =  \frac{2^{3-\frac{n+m}{2}}   \Gamma   \left(  \frac{n+m}{2}
  \right)}{\Gamma \left( \frac{n}{2} \right) \Gamma \left( \frac{m}{2}
  \right)    \,   \Gamma^2    \left(    \frac{m+n}{4}   \right)}    \:
r^{\frac{n+m}{2}-1} \, K^2_{\frac{n-m}{4}}(r)
\label{En_Student-t:eq}
\ee
where $K_\nu$ is  the modified Bessel function of  the third kind with
order $\nu$ (see for example \cite[8.432]{GraRyz80}).

Moreover, if $p > \frac{2 n}{n+m}$, the sum of the entropy rates is
\be\ba{lll}
\hspace{-.5cm}U_p(X_n) &  = & \log \pi  + \frac{(4 + q  \, (4-n-m)) \,
  \log 2}{2  \, n  \, (2-q)} +  \frac{2}{n \, (2-p)}  \left((p-1) \log
  \Gamma \left( \frac{n}{2} \right) \right.\vspace{2.5mm}
\\
&& \left.  +  \log \Gamma \left( \frac{(n+m) \, p -  2 n}{4} \right) -
  \log \Gamma  \left( \frac{(n + m)  \, p}{4} \right) +  p \log \Gamma
  \left( \frac{n+m}{4} \right) \right)\vspace{2.5mm}
\\
&&  +   \frac{2}{n  \,  (2-q)}   \displaystyle  \log  \int_0^{+\infty}
r^{\frac{(m-n) \, q}{4}+n-1} K_{\frac{n-m}{4}}^q(r) \, \d r
\label{U_Renyi_Student-t:eq}
\ea\ee
while, for $p=2$,
\be\ba{lll}
\hspace{-1cm} U_2(X_n) &  = & \log \pi + \frac{n-2}{n} \,  \log 2 + \,
\frac{2}{n}  \left(  \log \Gamma  \left(  \frac{m}{2}  \right) -  \log
  \Gamma   \left(   \frac{n+m}{2}  \right)   +   \log  \Gamma   \left(
    \frac{n+m}{4} \right) \right) \vspace{2.5mm}
\\
&& + \, \frac{n-m}{4  n} \left( \psi \left( \frac{n}{2}  \right) + 2 \,
  \psi  \left( \frac{n+m}{4}  \right)  \right) +  \frac{m}{n} \,  \psi
\left(  \frac{n+m}{2}  \right)  -  \frac{n+3m}{4  n}  \,  \psi  \left(
  \frac{m}{2} \right) \vspace{2.5mm}
\\
&& -  \,   \frac{2^{3-\frac{n+m}{2}}  \Gamma   \left(  \frac{n+m}{2}
  \right)}{n  \,  \Gamma  \left(  \frac{n}{2}  \right)  \Gamma  \left(
    \frac{m}{2}   \right)  \Gamma^2   \left(   \frac{n+m}{4}  \right)}
\displaystyle       \,       \int_0^{+\infty}      r^{\frac{n+m}{2}-1}
K^2_{\frac{n-m}{4}}(r) \, \log K^2_{\frac{n-m}{4}}(r) \, \d r
\label{U_Shannon_Student-t:eq}
\ea\ee

\end{theorem}

\begin{proof}
  (\ref{Dn_Student-t:eq})        is        straightforward        from
  (\ref{pdf_Student-t:eq}),  (\ref{fndn:eq})  and (\ref{Dn:eq})  while
  (\ref{En_Student-t:eq})  is  a   consequence  of  (\ref{En:eq})  and
  \cite[9-28~(22)]{Pou00}      (or      \cite[8.5~(20)]{Bat54}      or
  \cite[6.565--4]{GraRyz80}).             Finally,            plugging
  (\ref{Dn_Student-t:eq})     and     (\ref{En_Student-t:eq})     into
  (\ref{U_Renyi_si:eq})   and   using  \cite[8.380--3]{GraRyz80}   and
  $q=\frac{p}{p-1}$,     the    sum     of    the     entropy    rates
  (\ref{U_Renyi_Student-t:eq}) for $p \ne  2$ follows.  In the Shannon
  case,    starting    from    (\ref{U_Shannon_si:eq}),    recognizing
  $\displaystyle              \int_0^{+\infty}              \frac{r^{n
      -1}}{(1+r^2)^{\frac{n+m}{2}}} \,  \d r$  as a beta  integral and
  using    \cite[6.576--4]{GraRyz80}   to    evaluate   $\displaystyle
  \int_0^{+\infty} r^{-\lambda}  K_\nu^2(r) \,  \d r$, we  notice that
  $h(r)^\lambda   \log(h(r))   =   \frac{\partial}{\partial   \lambda}
  h(r)^\lambda$ (with  $h(r)=r$ and $h(r) = 1+r^2$)  to finally obtain
  (\ref{U_Shannon_Student-t:eq}).

  Notice  that  both  from  \cite[8.380--3]{GraRyz80}  to  insure  the
  existence  of  terms  in  (\ref{U_Renyi_Student-t:eq}) or  from  the
  asymptotics \cite[9.6.9 and 9.7.2]{AbrSte70} to insure the existence
  of the remaining integral, quantity $U_p(X_n)$ exists provided that
  \be p > \frac{2 n}{n+m} \label{Student-t_condition:eq} \ee
\end{proof}

Unfortunately,       in      the      general       case,      neither
(\ref{U_Renyi_Student-t:eq}) nor (\ref{U_Shannon_Student-t:eq}) can be
further analytically developed: numerical integration is necessary for
the  evaluation of  the remaining  integral.  We  will see  in section
\ref{Laplacian:sec} that  in the  special case $m=n+2$,  the remaining
integral can  be fully  developed and hence  the investigation  can be
completely performed.  However, an  asymptotic result can  be obtained
for any positive value of $m$, as expressed by the following theorem.

\begin{theorem}
\label{Student-t_asym:th}
For  any positive  value  of  the degree  of  freedom, the  following
equivalent holds:
\be
U_p(X_n) = \log(2  \pi) + \frac{\log p}{p-2} +  \frac{\log q}{q-2} +
o(1) \hspace{1cm} \mbox{for } p \ne 2
\label{Student-t_asym:eq}
\ee
\be
U_2(X_n) = 1 + \log \pi + o(1)
\label{Student-t_asym_2:eq}
\ee
Hence,   the   lower    bound   of   (\ref{uncertainty_Renyi:eq})   is
asymptotically attained when $n \to \infty$.

\end{theorem}

\begin{proof}
Let us denote
\be
I(\lambda,q)   =   \log    \left(   \int_0^{+\infty}   r^{-1}   \left(
    r^{\frac{m-n}{4}      +      \frac{n}{q}}     K_{\frac{n-m}{4}}(r)
  \right)^\lambda \, \d r \right)
\ee
Schwarz   inequality   implies   that  $\frac{\partial^2   I}{\partial
  \lambda^2} \ge  0$, showing  that function $I(\lambda,q)$  is convex
against $\lambda$.  As a consequence, for any $\lambda \in [1 ; 2]$ we
have $I(\lambda,q) = I((2-\lambda)  \times 1 + (\lambda-1) \times 2,q)
\le (2-\lambda) I(1,q)  + (\lambda-1) I(2,q)$.  If $\lambda  > 2$ this
inequality is reversed.  Since the remaining integral is $I(q,q)$, for
$q \ne 2$ we obtain the inequality
\be
\frac{1}{2-q}  \log \left( \int_0^{+\infty}  r^{\frac{(m-n) q}{4}+n-1}
  K^q_{\frac{n-m}{4}}(r)  \,   \d  r  \right)  \:  \le   \:  I(1,q)  +
\frac{q-1}{2-q} \, I(2,q)
\ee
Using \cite[6.561--16 and 6.576--4]{GraRyz80} to evaluate $I(1,q)$ and
$I(2,q)$ respectively, and  the fact that $p$ and  $q$ are conjugated,
we finally obtain, for $p \ne 2$,
\be\hspace{-5mm}\left\{\ba{ccl}
  U_p(X_n) & \le & M(n,m,p)\vspace{5mm}
  \\
  M(n,m,p) & = & \log(2 \pi)  + \frac{2}{n \, (2-p)} \left( (p-1) \log
    \Gamma \left( \frac{n}{2} \right) + \log \Gamma \left( \frac{(n+m)
        \, p - 2 n}{4} \right) \right.\vspace{2.5mm}
  \\
  & &  - \log \Gamma \left(  \frac{(n + m)  \, p}{4} \right) +  p \log
  \Gamma \left( \frac{n+m}{4} \right)\vspace{2.5mm}
  \\
  & & + (2-p) \, \log \Gamma \left( \frac{n (p-2) + m p}{4 p}\right) -
  2  \,  \log   \Gamma  \left(  \frac{2  n  (p-2)   +  (n+m)  p}{4  p}
  \right)\vspace{2.5mm}
  \\
  & & \left.  - \log \Gamma \left( \frac{n (p-2) + m p}{2 p} \right) +
    \log  \Gamma \left(  \frac{2 n  (p -  2) +  (n+m) p}{2  p} \right)
  \right)\vspace{2.5mm}
  \\
  & & + \frac{2}{n} \log \Gamma \left( \frac{n}{2 q} \right) + \frac{2
    (q-1)}{n (2-q)} \log \Gamma \left( \frac{n}{q} \right)
\ea\right.\ee
Now    using    the   asymptotics    of    the   log-gamma    function
\cite[6.1.41]{AbrSte70}  and tedious  algebra, one  obtains that the
upperbound $M$ verifies
\be
M(n,m,p) =  \log(2 \pi)  + \frac{\log p}{p-2}  + \frac{\log  q}{q-2} +
o(1)
\ee
whatever $m>0$, dependent  or not on $n$.  In the  case $p=2$, one has
by continuity (both for $U_p$ and for $M(n,m,p)$)
\be
U_2(X_n) \le M(n,m,2)
\ee
with
\bens\ba{lll}
\hspace{-5mm} M(n,m,2)  & =  & \log(2 \pi)  + \frac{2}{n}  \left( \log
  \Gamma \left(  \frac{n}{4} \right) - \log  \Gamma \left( \frac{n}{2}
  \right)  + \log  \Gamma  \left( \frac{m}{4}  \right)  - \log  \Gamma
  \left( \frac{n+m}{4} \right) \right)\vspace{2.5mm}
\\
& & + \frac{1}{2} \, \psi \left( \frac{n}{2} \right) -\frac{m}{2 \, n}
\, \psi \left( \frac{m}{2} \right) + \frac{m-n}{2 \, n} \, \psi \left(
  \frac{m+n}{2} \right) + \psi \left( \frac{m+n}{4} \right)
\ea\eens

Using again  the asymptotics of  the log-gamma and psi  functions, one
obtains
\be
M(n,m,2)   =  1 +  \log   \pi  +  o(1)
\ee
Together   with  (\ref{uncertainty_Renyi:eq})  these   results  confirm
(\ref{Student-t_asym:eq}):       the        lower       bound       of
(\ref{uncertainty_Renyi:eq}) is  attained when  $n \to \infty$  in the
Student-t case.
\end{proof}

%%%%%%%%%%%%%%%%%%%%%%%%%%%%%%%%%%%%%%%%%%%%%%%%%%%%%%%%%%%%%%%%%%%%%

\subsubsection{The Student-t case with $m = n+2$ degrees of freedom: a completely analytical study}
\label{Laplacian:sec}

Assuming $m=n+2$, the Bessel function in (\ref{En_Student-t:eq}) is of
order $- \, \frac{1}{2}$ and from \cite[8.469--3]{GraRyz80} we have,
\be
K_{- \frac{1}{2}}(r) = K_{ \frac{1}{2}}(r) = \sqrt{\frac{\pi}{2 r}} \,
\mbox{e}^{-r}
\ee

Notice first that  from the analogy with (\ref{fndn:eq})-(\ref{Dn:eq})
and using  the Euler's duplication  formula \cite[8.335--1]{GraRyz80},
$\widetilde{X}_n$ is distributed according to
\be
\widetilde{f}_n(x)  =  \frac{2^{n-1}   \,  \Gamma  \left(  \frac{n}{2}
  \right)}{\pi^{\frac{n}{2}} \Gamma(n)  } \, \exp  \left( - 2  \, (x^t
  x)^{\frac{1}{2}} \right).
\label{pdf_Laplacian:eq}
\ee
This  distribution is  called the  multivariate generalization  of the
Laplace                distribution\footnote{See               however
  \cite{FanKot90,KozPod00,KotKoz00}  for  a  different  definition  of
  generalized  multivariate Laplace  distributions.}   \cite{Ern98} or
multivariate   extension  generalized   Gaussian   \cite{CosHer03}  or
multivariate   exponential  power   \cite{KotKoz00}.    As  shown   in
\cite{MatKri96}, each particle of  an ideal relativistic photon gas in
a  container  with rigid  and  diathermic  walls  has a  3-dimensional
momentum X that follows distribution (\ref{pdf_Laplacian:eq}).

\begin{theorem}
  The  entropy sum  in the  case  of the  Student-t distribution  with
  $m=n+2$ degrees of freedom is, for $p \ne 2$,
\be\ba{lll}
\hspace{-6mm}\displaystyle  U_p(X_n) &  = &  \displaystyle \log  \pi +
\frac{2 \log  q - q \log 2}{q-2}  + \frac{2}{n (2-p)} \left(  \log 2 +
  \log \Gamma(n) \right.\vspace{2.5mm}
\\
&&  \left.  -  \log \Gamma  \left( \frac{n}{2}  \right) +  \log \Gamma
  \left(  \frac{(n+1)  \, p  -  n}{2}  \right)  - \log  \Gamma  \left(
    \frac{(n+1) \, p}{2} \right) \right)
\label{U_Renyi_Laplacian:eq}
\ea\ee
and, for $p=2$,
\be
U_2(X_n) = 1 + \log \left(\frac{\pi}{2} \right) + \frac{n+1}{n} \left(
  \psi(n) - \psi \left( \frac{n}{2} \right) - \frac{1}{n} \right)
\label{U_Shannon_Laplacian:eq}
\ee
Moreover, it holds asymptotically for $p \ne 2$
\be
\frac{H_{\frac{p}{2}}(X_n)  + H_{\frac{q}{2}}(\widetilde{X}_n)  }{n} =
\log (2 \pi) + \frac{\log p}{p-2} + \frac{\log q}{q-2} + o(1)
\label{U_Renyi_Laplacian_asym:eq}
\ee
while for $p=2$,
\be
\frac{H(X_n) + H(\widetilde{X}_n)}{n} = 1 + \log \pi + o(1)
\label{U_Shannon_Laplacian_asym:eq}
\ee

\end{theorem}

\begin{proof}
  (\ref{U_Renyi_Laplacian:eq})    can   be    easily    deduced   from
  (\ref{U_Renyi_Student-t:eq}) using the fact that $q = \frac{p}{p-1}$
  and  Euler's  duplication  formula  \cite[8.335--1]{GraRyz80}.   The
  limit   when    $p   \rightarrow    2$   can   be    obtained   from
  (\ref{U_Shannon_Student-t:eq}) or  using a first  order expansion of
  (\ref{U_Renyi_Laplacian:eq})  and formula \cite[8.365--1]{GraRyz80}.
  This last  result can also  be found starting directly  from Shannon
  entropy.      Using      \cite[6.1.41     and     6.3.18]{AbrSte70},
  (\ref{U_Renyi_Laplacian_asym:eq})                                 and
  (\ref{U_Shannon_Laplacian_asym:eq}) follow straightforwardly.
\end{proof}

These results  show that the Student-t distribution  with $m=n+2$ (or,
by    conjugation,   the    multivariate   exponential    power   pdf
(\ref{pdf_Laplacian:eq}))   achieves    asymptotically   equality   in
(\ref{uncertainty_Renyi:eq})   (and   (\ref{uncertainty_Shannon:eq})).
Contrarily  to  the  result  of  Abe \cite{AbeRaj01}  for  the  Cauchy
distribution  and Shannon entropy,  the results  obtained here  are in
analytical  form.   We will  discuss  in section  \ref{discussion:sec}
about the possible explanations of this asymptotic behavior.

%%%%%%%%%%%%%%%%%%%%%%%%%%%%%%%%%%%%%%%%%%%%%%%%%%%%%%%%%%%%%%%%%%%%%

\subsubsection{The Student-t case with 1 degree of freedom: revisiting the Cauchy case}
\label{Cauchy:sec}

We  deal  in this  section  with the  case  of  a multivariate  Cauchy
distribution as previously studied by Abe in \cite{AbeRaj01}, \ie
\be
f_n(x)        =        \frac{\Gamma        \left(        \frac{n+1}{2}
  \right)}{\pi^{\frac{n+1}{2}} }( 1 + x^t x)^{-\frac{n+1}{2}}
\label{pdf_Cauchy:eq}
\ee

Formulas                (\ref{U_Renyi_Student-t:eq})               and
(\ref{U_Shannon_Student-t:eq})  apply with  $m=1$:  unfortunately, the
remaining  integrals in  this case  cannot be  further  simplified and
numerical evaluations  are required. Notice that in  the Shannon case,
the result found  by Abe in \cite{AbeRaj01} is  recovered, except that
one of the two integrals, numerically evaluated in \cite{AbeRaj01}, is
in fact analytically expressed here.

\begin{figure}[htbp]
\psfrag{U2}{\footnotesize $U_2(X_n)$}
\psfrag{U3}{\footnotesize $U_3(X_n)$}
\psfrag{U10}{\footnotesize $U_{10}(X_n)$}
\psfrag{n}{\footnotesize $n$}
\includegraphics[width=4.5cm]{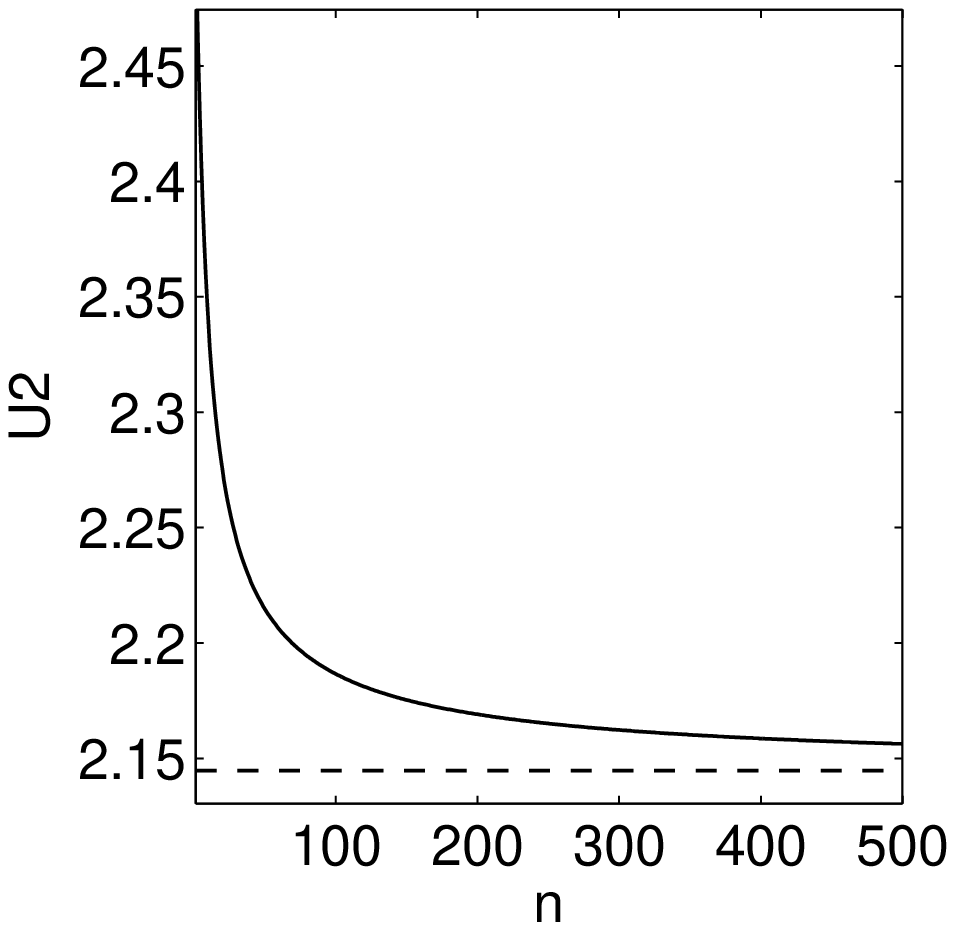}
\includegraphics[width=4.5cm]{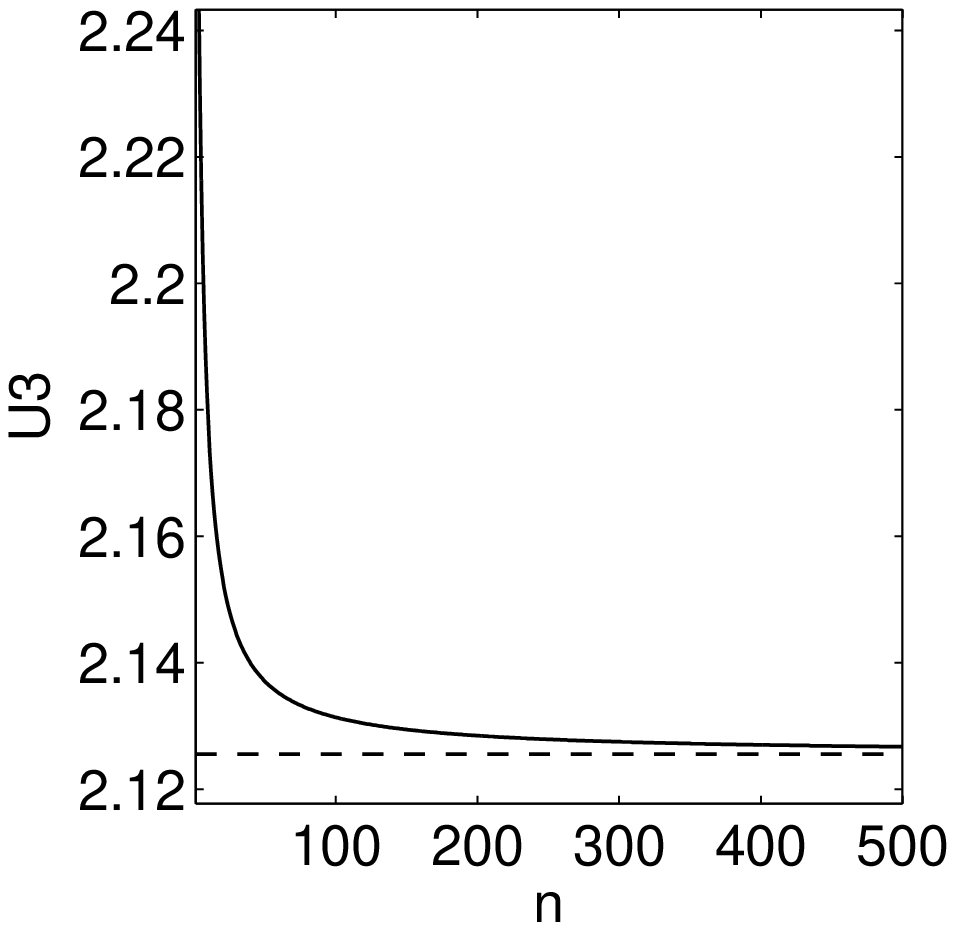}
\includegraphics[width=4.5cm]{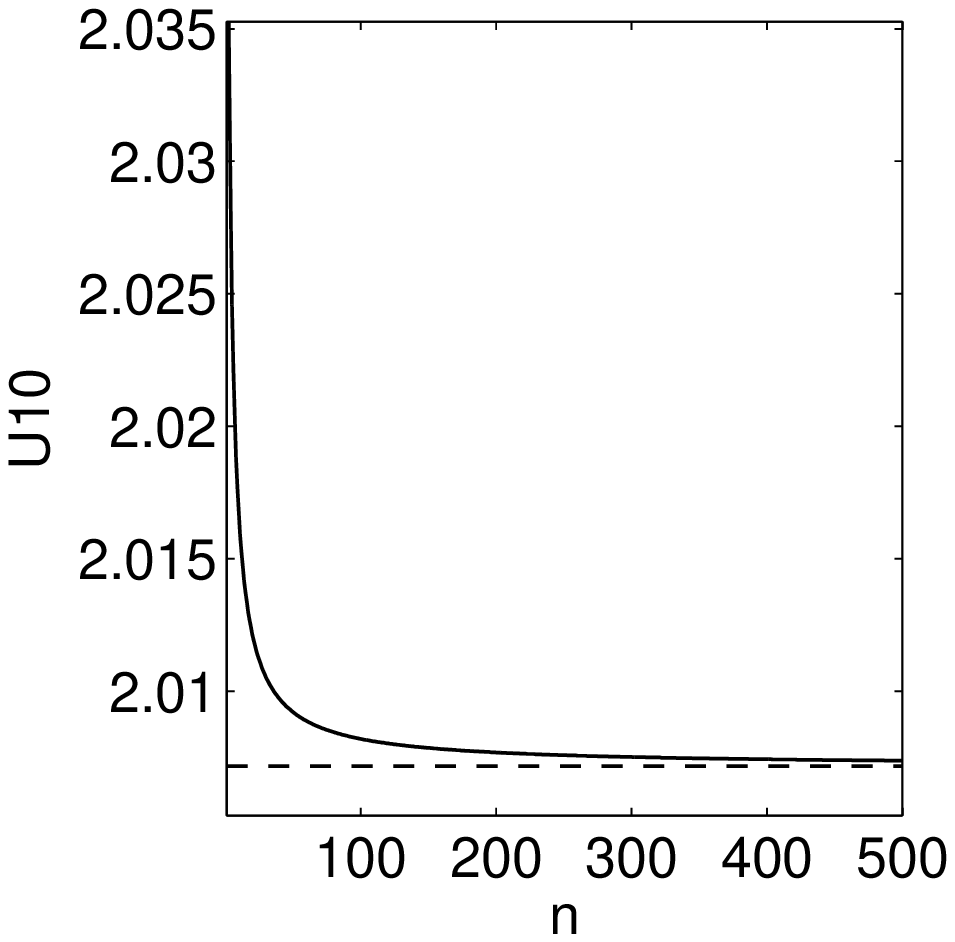}
\caption{Illustration      of       the      uncertainty      relation
  (\ref{uncertainty_Renyi:eq})  in  the   Cauchy  context,  for  $p=2$
  (Shannon  version  (\ref{uncertainty_Shannon:eq})),  $p=3$  and  for
  $p=10$ respectively.  The solid line  depicts the sum of the entropy
  rates  $U_p(X_n)$  as  a  function  of  $n$,  and  the  dashed  line
  represents the lower bound.}
\label{Cauchy_Renyi:fig}
\end{figure}

Curves  in   Fig.   \ref{Cauchy_Renyi:fig}  depict   the  behavior  of
(\ref{U_Shannon_Student-t:eq})  and (\ref{U_Renyi_Student-t:eq})  as a
function  of $n$,  in the  Cauchy case,  (the integral  is numerically
evaluated), for  $p=2$ (Shannon case), $p=3$  and $p=10$ respectively.
The same behavior is observed for  any value of $p$ tested.  The first
curve  simply confirms the  results by  Abe \cite{AbeRaj01}.   It also
appears from  these curves  that the conclusion  found by  Abe remains
true for the R\'enyi version of the uncertainty relation, whatever the
possible  values of  $p$ and  thus  confirm (\ref{Student-t_asym:eq}).
Hence, this behavior is not specific to the Shannon entropy.

%%%%%%%%%%%%%%%%%%%%%%%%%%%%%%%%%%%%%%%%%%%%%%%%%%%%%%%%%%%%%%%%%%%%%

\subsection{The general Student-r case}
\label{Student-r:sec}

We consider  now the  case where $X_n$  is distributed according  to a
Student-r law with $m$  degree of freedom\footnote{$m$ is named degree
  of  freedom  by   misuse  of  language  and  may   be  called  shape
  parameter.  This choice is adopted by analogy  with  the Student-t
  case.
}, \ie

\be
f_n(x) = \frac{\Gamma \left( \frac{m}{2} +1 \right)}{\pi^{\frac{n}{2}}
  \Gamma  \left(  \frac{m-n}{2} +  1  \right)} \,  \left(  1  - x^t  x
\right)_+^{\frac{m-n}{2}}
\label{pdf_Student-r:eq}
\ee
where $ m > n-2$ and $(.)_+ = \max(.,0)$.  When $m=n$, vector $X_n$ is
uniformly distributed in the $n$-dimensional  sphere $x^t x = 1$.  The
Student-r variables play also  an important role in probability, since
they appear  as $n-$dimensional marginals of  the uniform distribution
on the sphere in $\Rset^{m+2}$; they are also maximizers (for $\lambda
= 1 + \frac{2}{m-n}$)  of the R\'enyi/Tsallis entropy under covariance
constraint \cite{ParBir05,TsaMen98,Tsa99,VigHer04}.

\begin{theorem}
  If $X_n$ is Student-r distributed, then the pdf of $\| X_n \|$ is
\be
D_n(r) =  \frac{2 \Gamma \left( \frac{m}{2}  +1 \right)}{\Gamma \left(
    \frac{n}{2}  \right) \Gamma  \left( \frac{m-n}{2}  +1  \right)} \,
r^{n-1} (1-r^2)^{\frac{m-n}{2}}
\label{Dn_Student-r:eq}
\ee
for  $r  \in  [0 ;  1]$  and  $0$  otherwise,  while  the pdf  of  $\|
\widetilde{X}_n \|$ is
\be
E_n(r)  =  \frac{2^{\frac{m-n}{2}+1}   \Gamma  \left(  \frac{m}{2}  +1
  \right)  \Gamma^2  \left(  \frac{m-n}{4} +1  \right)}{\Gamma  \left(
    \frac{n}{2}  \right) \Gamma  \left( \frac{m-n}{2}  +1  \right)} \,
r^{- \, \frac{m-n}{2}-1} J^2_{\frac{m+n}{4}}(r)
\label{En_Student-r:eq}
\ee
for $r \in [0 ; +\infty)$  and $0$ otherwise.

Moreover, if $q > \frac{4 n}{m+n+2}$, the sum of the entropy rates is
\be\ba{lll}
U_p(X_n) & = & \displaystyle \log \pi +  \frac{(4 + ( m - n ) \, q) \,
  \log 2}{2 \, n \, (2-q)}\vspace{2.5mm}
\\
&  &  +  \frac{2}{n \,  (2-p)}  \left(  (p-1)  \, \log  \Gamma  \left(
    \frac{n}{2} \right)  + \log \Gamma \left( \frac{(m-n)  \, p}{4} +1
  \right) \right.\vspace{2.5mm}
\\
&& \left.   - \log \Gamma \left(  \frac{(m - n) \,  p + 2  \, n}{4} +1
  \right)   -  p   \log   \Gamma  \left(   \frac{m-n}{4}  +1   \right)
\right)\vspace{2.5mm}
\\
&&  +   \frac{2}{n  \,  (2-q)}   \displaystyle  \log  \int_0^{+\infty}
r^{\frac{- (m+n)  \, q}{4}+n-1} \left|  J_{\frac{m+n}{4}}(r) \right|^q
\, \d r
\label{U_Renyi_Student-r:eq}
\ea\ee
for any $m > n-2$, while, for $p=2$,
\be\ba{lll}
%
%\lefteqn{
\hspace{-5mm}U_2(X_n) & = & \log  (2 \pi) + \, \frac{2}{n} \left( \log
  \Gamma  \left(  \frac{m-n}{2}  +  1  \right) -  \log  \Gamma  \left(
    \frac{m}{2} + 1 \right) \right.\vspace{2.5mm}
\\
&& \left. -  \log \Gamma \left( \frac{m-n}{4} +1  \right) \right) + \,
\frac{m+n}{4 n}  \left( \psi  \left( \frac{n}{2} \right)  + 2  \, \psi
  \left( \frac{m-n}{4} +1 \right) \right)\vspace{2.5mm}
\\
&&  -   \frac{m}{n}  \,  \psi   \left(  \frac{m-n}{2}  +1   \right)  +
\frac{3m-n}{4     n}     \,     \psi     \left(     \frac{m}{2}     +1
\right)\nonumber\vspace{2.5mm}
\\
&& - \, \frac{2^{\frac{m-n}{2}+1} \Gamma \left( \frac{m}{2} +1 \right)
  \,  \Gamma^2 \left(  \frac{m-n}{4}  +1 \right)}{n  \, \Gamma  \left(
    \frac{n}{2}  \right)  \Gamma   \left(  \frac{m-n}{2}  +1  \right)}
\displaystyle           \int_0^{+\infty}          r^{-\frac{m-n}{2}-1}
J^2_{\frac{m+n}{4}}(r) \, \log J^2_{\frac{m+n}{4}}(r) \, \d r
\label{U_Shannon_Student-r:eq}
\ea\ee
\end{theorem}

\begin{proof}
  (\ref{Dn_Student-r:eq})      is       directly      issued      from
  (\ref{pdf_Student-r:eq})  while (\ref{En_Student-r:eq})  is obtained
  from      (\ref{En:eq})      and     \cite[9-28~(3)]{Pou00}      (or
  \cite[8.5~(33)]{Bat54}    or    \cite[6.567--1]{GraRyz80}).    Then,
  plugging        expressions        (\ref{Dn_Student-r:eq})       and
  (\ref{En_Student-r:eq})   into   (\ref{U_Renyi_si:eq})   and   using
  \cite[8.380--1]{GraRyz80} yields (\ref{U_Renyi_Student-r:eq}) for $p
  \ne 2$ and for any $m > n-2$.

  In    the    Shannon    case,    we    can    again    start    from
  (\ref{U_Shannon_si:eq}).   Using the  same technique  as  in section
  \ref{Student-t:sec}      and      result     \cite[8.380--1      and
  6.574--2]{GraRyz80}, we obtain (\ref{U_Shannon_Student-r:eq}).

  From the asymptotics \cite[9.1.7  and 9.2.1]{AbrSte70} of the Bessel
  function  for  small  and  large  argument,  the  integral  converges
  provided that
  \be q > \frac{4 n}{m+n+2} \label{Student-r_condition:eq}\ee
\end{proof}

Again, in  the general case,  neither (\ref{U_Renyi_Student-r:eq}) nor
(\ref{U_Shannon_Student-r:eq})  can be further  analytically developed
and recourse  to numerical integration  for the remaining  integral is
needed.           Note          the          symmetry          between
(\ref{U_Renyi_Student-r:eq})-(\ref{U_Shannon_Student-r:eq})         and
(\ref{U_Renyi_Student-t:eq})-(\ref{U_Shannon_Student-t:eq})
respectively that  can be explained by remarking  the symmetry between
the Student-t and Student-r variables as evoked in \cite{CosHer03}.

%%%%%%%%%%%%%%%%%%%%%%%%%%%%%%%%%%%%%%%%%%%%%%%%%%%%%%%%%%%%%%%%%%%%%

\subsubsection{The Student-r case with $m = n$ degrees of freedom: the uniform case}
\label{uniform:sec}

\begin{figure}[htbp]
\psfrag{U2}{\footnotesize $U_{1.91}(X_n)$}
\psfrag{U3}{\footnotesize $U_{1.5}(X_n)$}
\psfrag{U10}{\footnotesize $U_{1.11}(X_n)$}
\psfrag{n}{\footnotesize $n$}
 \includegraphics[width=4.5cm]{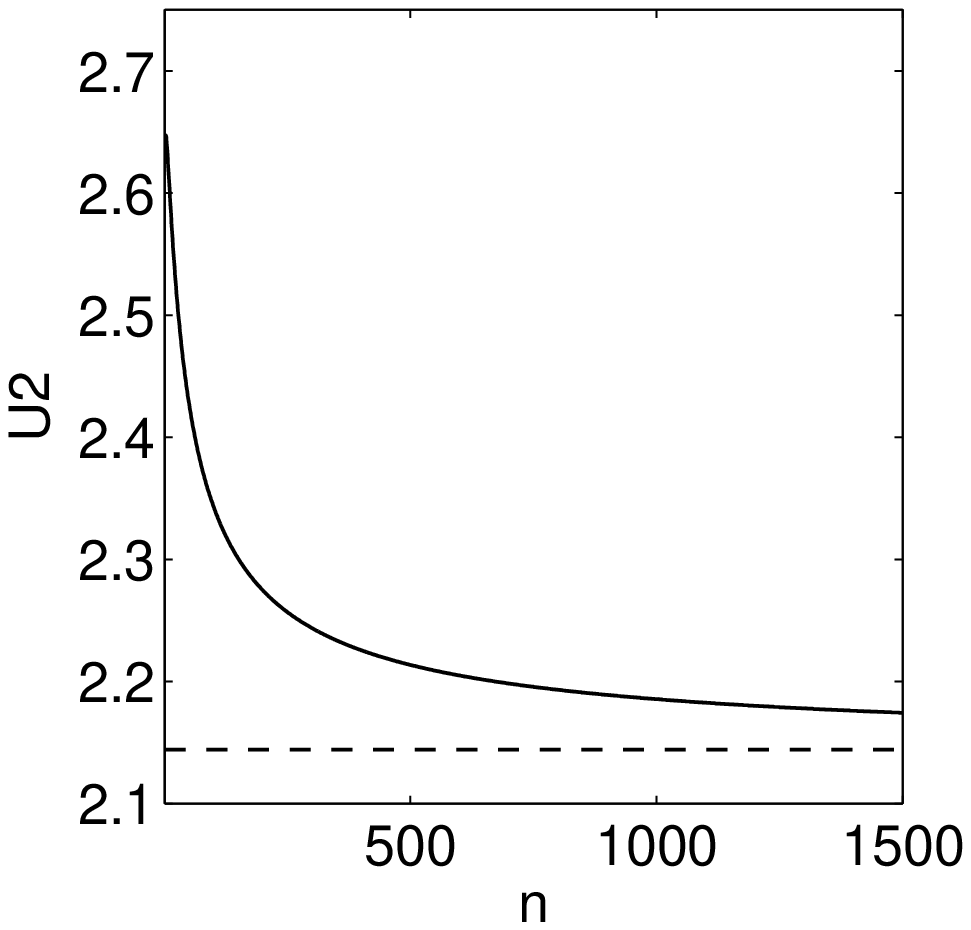}
 \includegraphics[width=4.5cm]{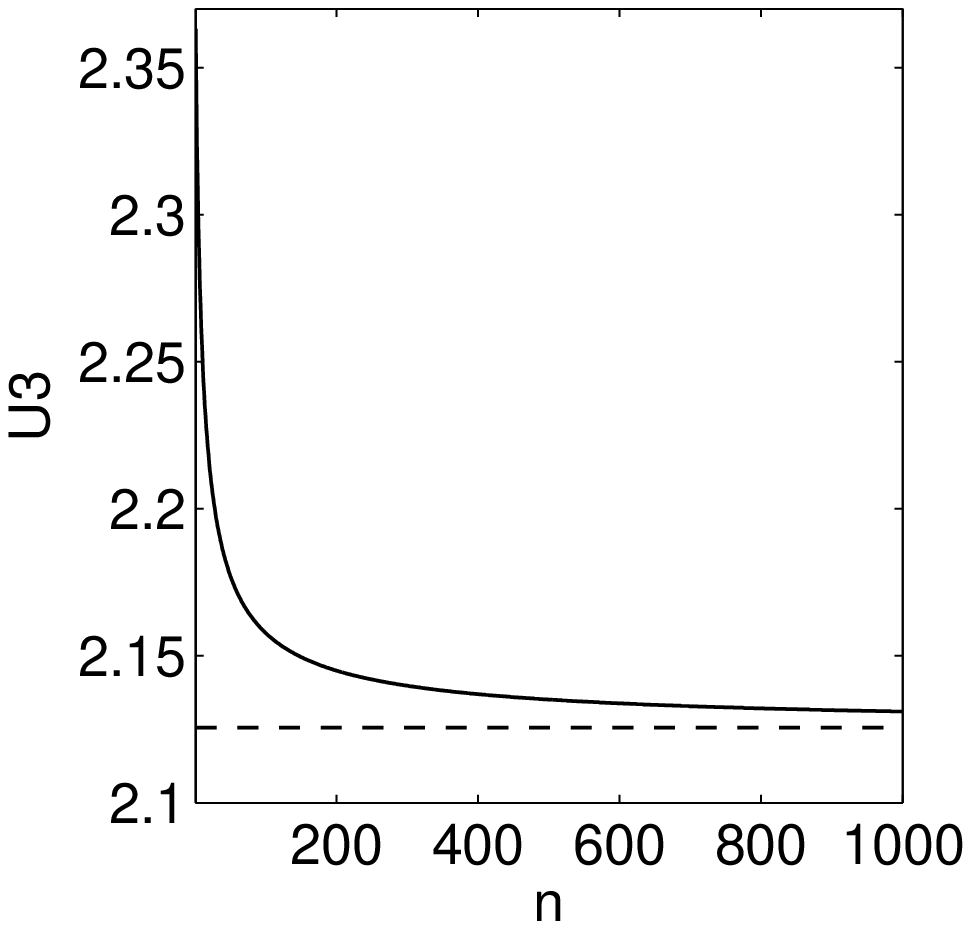}
 \includegraphics[width=4.5cm]{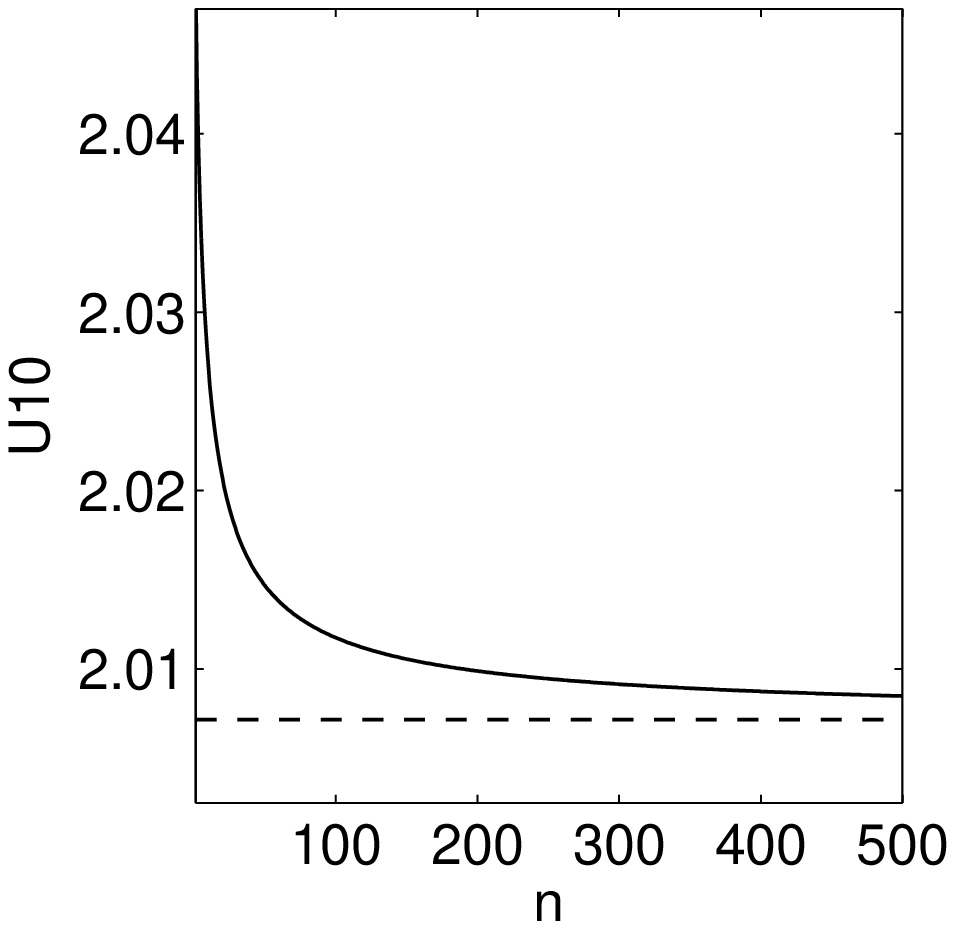}
 \caption{Illustration      of      the      uncertainty      relation
   (\ref{uncertainty_Renyi:eq}) in the uniform case, for $q=2.1$ (near
   the Shannon version  (\ref{uncertainty_Shannon:eq})), $q=3$ and for
   $q=10$  respectively.   The  solid  line depicts  $U_p(X_n)$  as  a
   function of $n$ and the dashed line depicts the lower bound.}
\label{Uniform_Renyi:fig}
\end{figure}

In this case,  we directly apply formulas (\ref{U_Renyi_Student-r:eq})
and (\ref{U_Shannon_Student-r:eq}), where  the remaining integrals are
numerically  evaluated.  Figures \ref{Uniform_Renyi:fig}  depicts then
the      behavior      of      (\ref{U_Shannon_Student-r:eq})      and
(\ref{U_Renyi_Student-r:eq}) as  a function of $n$,  for $q=2.1$ (near
the Shannon  case), $q=3$ and $q=10$ respectively.   The same behavior
is observed for any value of  $p$ tested, showing again that the lower
bounds  of the uncertainty  relations are  attained as  $n$ increases.
Again, we  currently try to  end analytically the  investigation. Note
that we  observe an increasingly slower convergence  as $q$ approaches
$2$ (especially for ``small'' $m$,  \eg uniform): we choose to present
the case  $q=2.1$ since  the convergence is  not too slow  compared to
other values of $q$.

%%%%%%%%%%%%%%%%%%%%%%%%%%%%%%%%%%%%%%%%%%%%%%%%%%%%%%%%%%%%%%%%%%%%%

\subsubsection{Some other  Student-r cases}
\label{other_Student-r:sec}

\begin{figure}[htbp]
\psfrag{U2}{\footnotesize $U_{1.91}(X_n)$}
\psfrag{U3}{\footnotesize $U_{1.5}(X_n)$}
\psfrag{U10}{\footnotesize $U_{1.11}(X_n)$}
\psfrag{n}{\footnotesize $n$}
\psfrag{mm2}{\tiny $m = n+2$}
\psfrag{m3n}{\tiny $m = 2 n$}
\includegraphics[width=4.5cm]{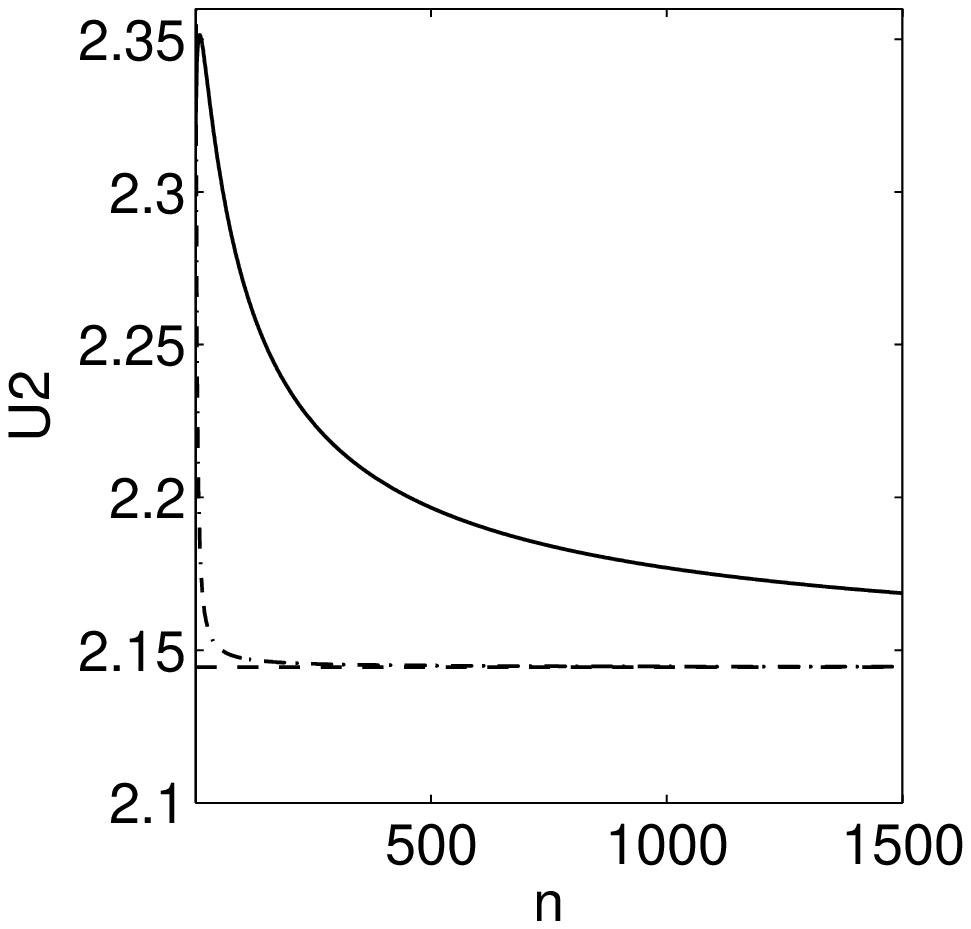}
\includegraphics[width=4.5cm]{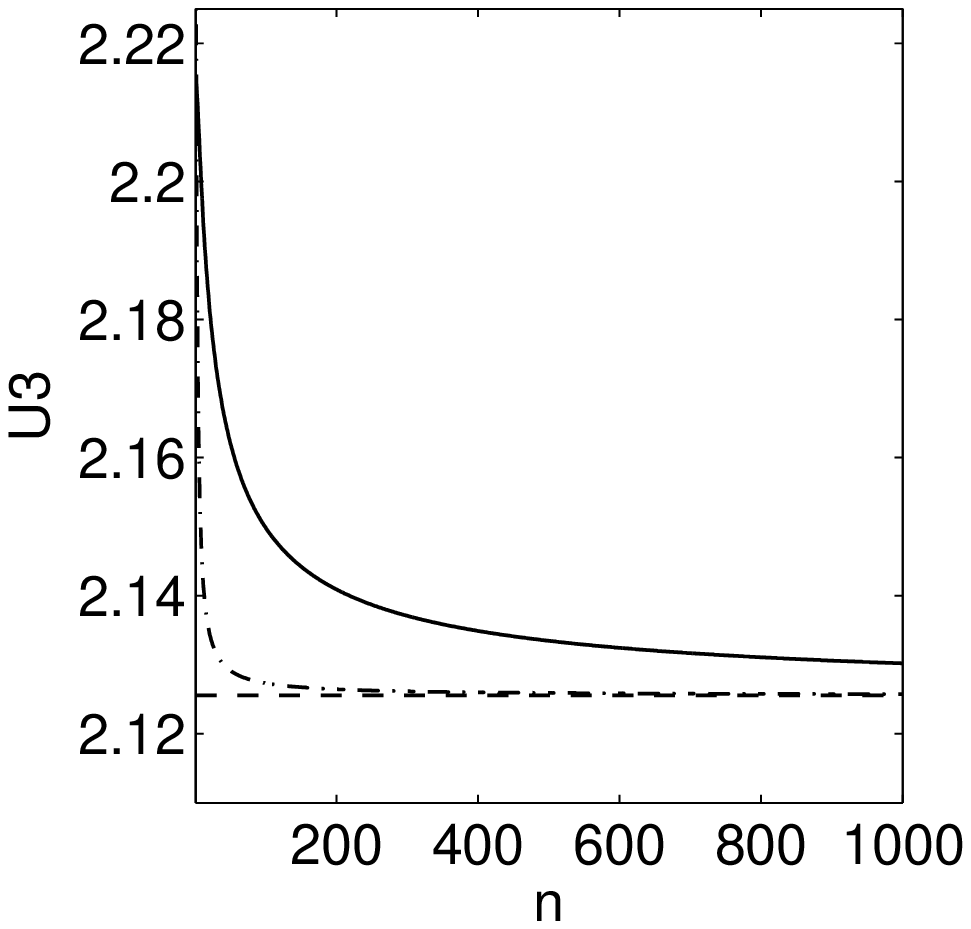}
\includegraphics[width=4.5cm]{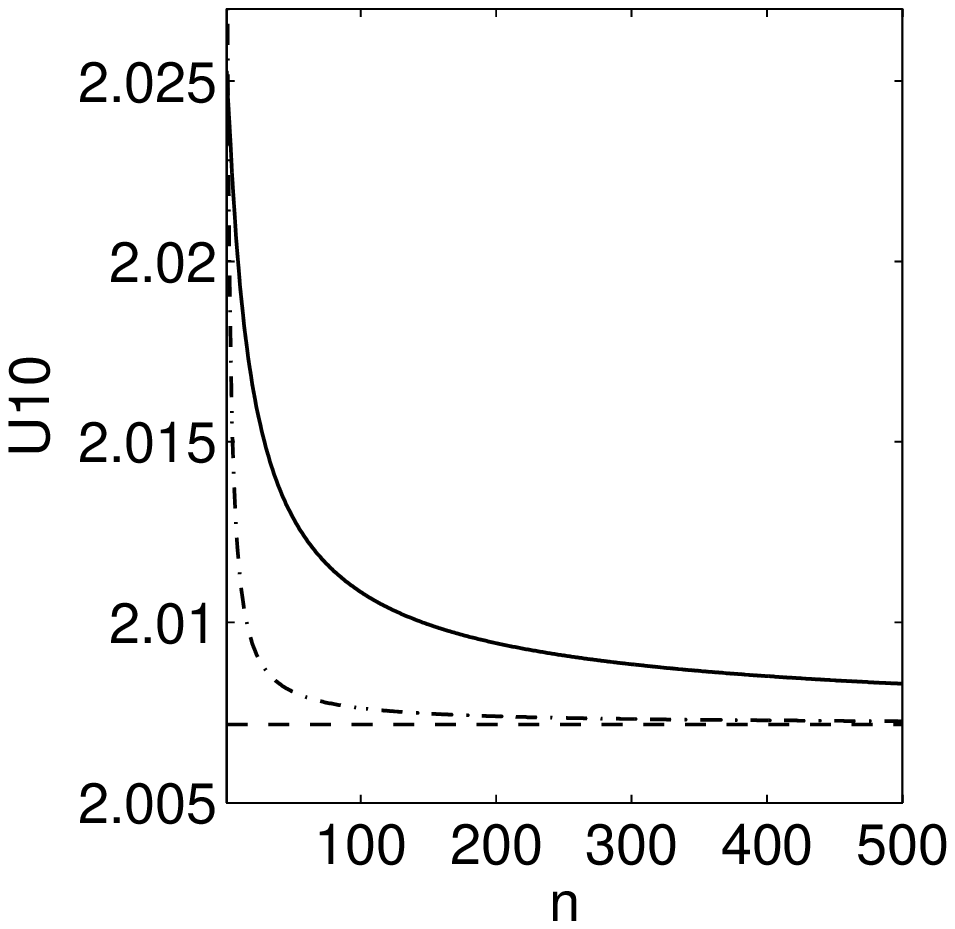}
\caption{Illustration      of       the      uncertainty      relation
  (\ref{uncertainty_Renyi:eq})  in  several  Student-r  contexts,  for
  $q=2.1$  (near the Shannon  version (\ref{uncertainty_Shannon:eq})),
  $q=3$ and for $q=10$  respectively.  The figures represent the cases
  $m=n+2$ (solid line) and $m=2 \, n$ (dashed-dotted line).  The small
  dashed line represents the lower bound.}
\label{Student-r_Renyi:fig}
\end{figure}

Again we do not enter  into details in this subsection. However, other
cases were  studied as illustrated  in Fgure \ref{Student-r_Renyi:fig}
for $m = n + 2$ and $m = 2 n$ and parameters $q = 2$, $q = 3$ and $q =
10$ respectively, showing that the asymptotic behavior holds. Contrary
to  the Student-t  case, in  the  Student-r case  we do  not find  any
specific  case where a  completely analytical  study can  be performed
although $J_{k+\frac{1}{2}}$  admits explicit formulation  when $k \in
\Nset$ \cite[8.462]{GraRyz80}.

%%%%%%%%%%%%%%%%%%%%%%%%%%%%%%%%%%%%%%%%%%%%%%%%%%%%%%%%%%%%%%%%%%%%%

\subsection{Discussion}
\label{discussion:sec}

The  preceding results show  that any  Student-t or  Student-r vectors
reach asymptotically the case  of equality in the Bialynicki-Birula \&
Mycielski inequality or  in its R\'enyi extension.  As  the only exact
(finite dimensional) case of equality  is met by Gaussian vectors, one
may   be  tempted   to   explain  this   asymptotic   behavior  by   a
``Gaussianization effect'' of these vectors, namely the fact that they
become ``more and  more Gaussian'' in some sense  as $n$ increases. In
the  rest  of  this  paragraph,  we study  two  possible  measures  of
Gaussianization:  in  the   distribution  sense  and  the  information
divergence sense.

%%%%%%%%%%%%%%%%%%%%%%%%%%%%%%%%%%%%%%%%%%%%%%%%%%%%%%%%%%%%%%%%%%%%%

\subsubsection{Convergence in distribution}
\label{distribution:sec}

A well-known  property of a Student-t  vector $X_n$ is that  it can be
expressed    as    a     Gaussian    scale    mixture    \cite{Chu73},
namely\footnote{$\egald$ means equality in distribution.}
\be
X_n \egald \sqrt{A_m} \times G_n,
\label{Student-t_representation:eq}
\ee
where  $A_m$ is  a scalar  inverse Gamma  random variable\footnote{The
  inverse  Gamma distribution  $\mbox{inv}\Gamma(\alpha,\beta)$ writes
  $f(x)  = \frac{\beta^{\alpha}  }{\Gamma(\alpha)}  x^{-1-\alpha} \exp
  \left(-\frac{\beta  }{x}\right)$  it  is  the  distribution  of  the
  inverse    of   a    $\Gamma(\alpha,1/\beta)$    random   variable.}
$\mbox{inv}\Gamma \left( \frac{m}{2} , 2 \right)$ with shape parameter
$m/2$ and scale parameter $1/2$, independent of the zero-mean Gaussian
vector   $G_n$    with   identity   covariance    matrix   (see   also
\cite{VigHer04,CosHer03,VigBer03}).   In  the  particular Cauchy  case
$m=1$,  one  recovers  the  fact  that  $A_1$  is  a  L\'evy  variable
\cite{SamTaq94}.  Since  $E \left[ \sqrt{m-2} \,  \sqrt{A_m} \right] =
\sqrt              {\frac{m-2}{2}}              \,              \frac{
  \Gamma(\frac{m-1}{2})}{\Gamma(\frac{m}{2})} =  1 + O(1/m)$  and $VAR
\left[ \sqrt{m-2} \, \sqrt{A_m} \right]  = 1 - \frac{m-2}{2} \, \left(
  \frac{\Gamma \left( \frac{m-1}{2} \right)}{\Gamma \left( \frac{m}{2}
    \right)}  \right)^2  =  \frac{1}{2m}  + o(1/m)$,  we  deduce  that
$\sqrt{m-2}  \, \sqrt{A_m}$ converges  to $1$  almost surely  when $m$
goes  to infinity.   Hence, by  Slutsky's  theorem, when  $m$ goes  to
infinity  with  $n$,  any  subvector  $X^k_n$  of  $X_n$  with  finite
dimension  $k<n$ converges  in distribution  to a  Gaussian  vector of
finite size $k$ \cite{Muk00}.

However, when $m$ remains constant and $n$ grows, since variable $A_m$
does not depend  on the dimension $n$, any  subvector $X_n^k$ of $X_n$
of  finite dimension  $k <  n$ remains  Student-t with  $m$  degree of
freedom: as $n \rightarrow \infty, $ no Gaussianization effect happens
for constant $m$, at least  in the distribution sense.  As an example,
this is the  case for a Cauchy vector for which  it is well-known that
any subvector $X_n^k$ remains Cauchy distributed whatever $n$.

\begin{figure}[htbp]
\psfrag{x}{\footnotesize $x$}
\psfrag{x1}{\footnotesize $x_1$}
\psfrag{x2}{\footnotesize $x_2$}
\psfrag{fn1}{\footnotesize $f_{n,1}(x)$}
\psfrag{fn2}{\footnotesize $f_{n,1,2}(x_1,x_2)$}
\includegraphics[width=5cm]{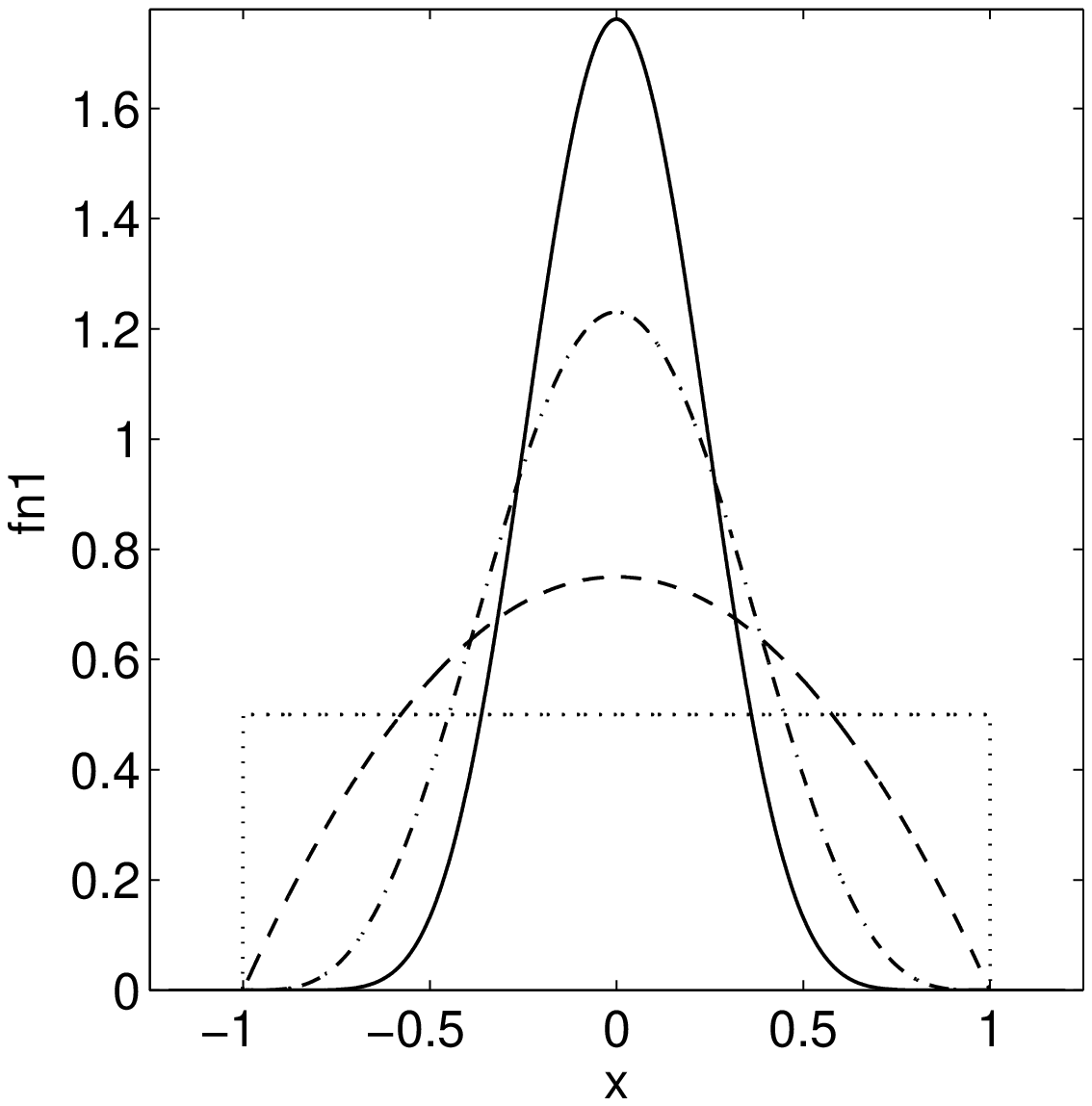}\hspace{1cm}
\includegraphics[width=6.5cm]{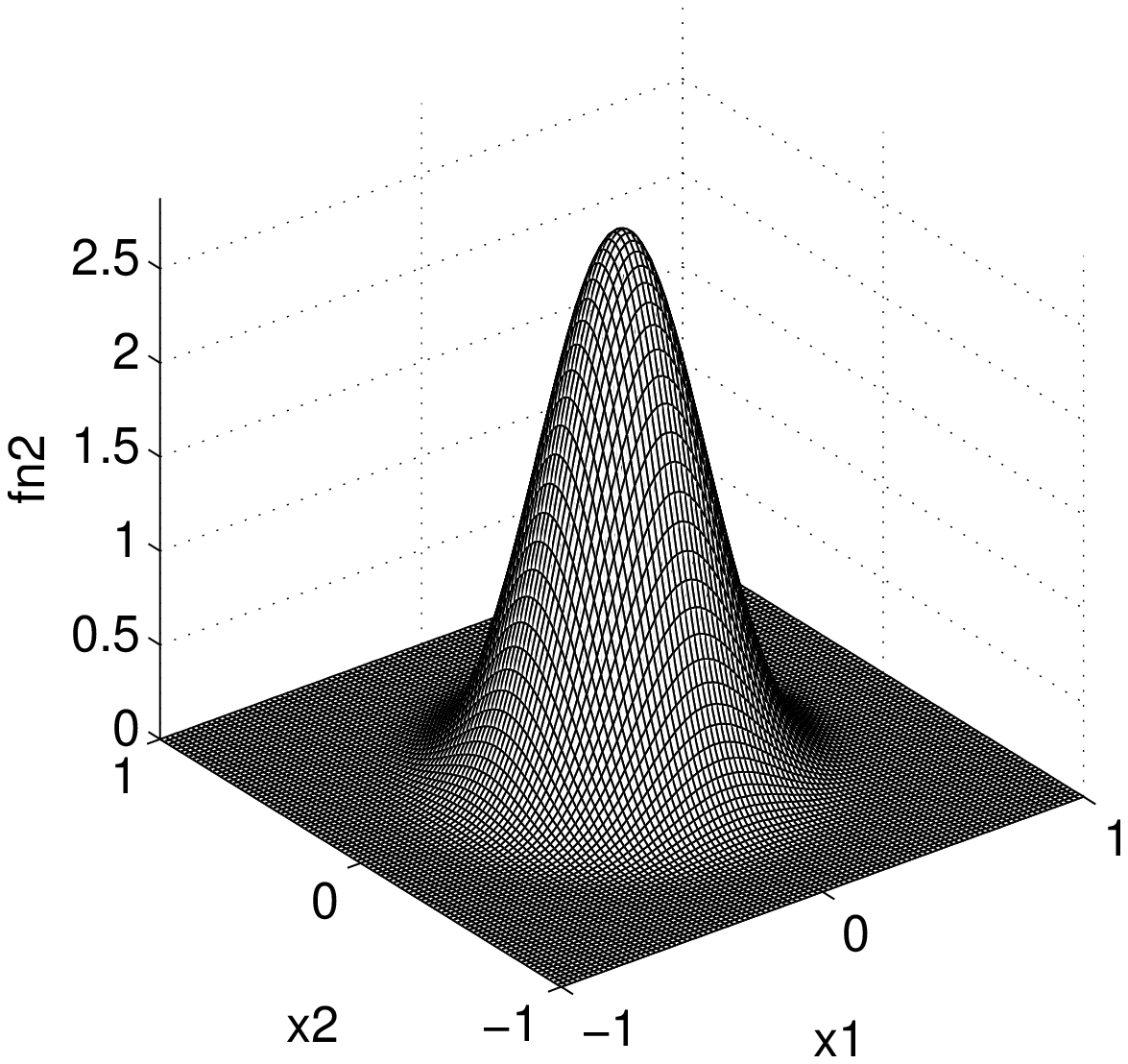}
\caption{Illustration  of the  convergence  in distribution  of the  a
  finite size  uniformly distributed vector  to the Gaussian  pdf: pdf
  $f_{n,1}(x)$ of  the first component  and for several values  of $n$
  (left) and pdf $f_{n,1,2}(x_1,x_2)$  of the first two components for
  $n=10$  (right).  The  dotted line  depicts $n=1$,  the  dashed line
  represent $n=2$,  the dashed line  depicts $n=5$ and the  solid line
  represents the case $n=10$.}
\label{Uniform_Gaussian:fig}
\end{figure}

For  the Student-r  random vectors,  the scale  mixture representation
does  not hold.   However, it  is  shown in  \cite{CosHer03} that  any
Student-r vector $X_n$ can be expressed as
\be
X_n    \egald    \frac{G_n}{\left(    \|    G_n   \|^2    +    B_{m,n}
  \right)^{\frac{1}{2}}}
\label{Student-r_representation:eq}
\ee
where  $B_{m,n}$ is a  scalar Gamma  distributed variable\footnote{The
  Gamma    distribution   $\Gamma(\alpha,\beta)$   writes    $f(x)   =
  \frac{1}{\beta^{\alpha} \Gamma(\alpha)} x^{\alpha-1} \exp \left(- \,
    \frac{x}{\beta}  \right)$.},   with  shape  parameter   $\alpha  =
(m-n+2)/2$ and  scale parameter $\beta  = 2$, which is  independent on
the unit  covariance Gaussian  vector $G_n$.  This  representation was
given in \cite{CosHer03} for $m$ integer, but it holds for non integer
values of $m$ as well.  Now, variable $C_m = \| G_n \|^2 + B_{m,n}$ is
Gamma distributed  $\Gamma\left( \frac{m+2}{2} , 2  \right)$ and then,
with    the   same    technique    as   in    the   Student-t    case,
$\sqrt{\frac{m+2}{C_m}}\to 1$ almost surely when $n \to \infty$ (since
$m > n-2$, $m+2 \to \infty$  when $n \to \infty$).  Although $C_m$ and
$G_n$  are not  independent,  one can  again  evoke Slutsky's  theorem
\cite{Muk00} to  conclude that a  finite-size subvector of  $X_n^k$ of
$X_n$ tends  in distribution  to a Gaussian  vector when $n$  tends to
infinity.  This  result is  known as ``Poincar\'e's  observation'' and
gave  birth  to an  important  literature;  despite  its name,  it  is
attributed   to   Borel   in   \cite{DiaFre87}  and   to   Mehler   in
\cite{Joh03}. It  is illustrated in  figure \ref{Uniform_Gaussian:fig}
in the uniform case for $k=1$ and various values of $n$, and for $k=2$
and $n=10$.

%%%%%%%%%%%%%%%%%%%%%%%%%%%%%%%%%%%%%%%%%%%%%%%%%%%%%%%%%%%%%%%%%%%%%

\subsubsection{Convergence in the Kullback-Leibler divergence rate sense}
\label{KL_divergence:sec}

The Kullback-Leibler  (KL) divergence between a  random vector $Y_n$
and a random  vector $Z_n$ with respective pdfs  $\rho_y$ and $\rho_z$
is defined as
\be
D_\kl(Y_n \,  \| \, Z_n)  = D_\kl(\rho_y \,  \| \,
\rho_z) = \int \rho_y \log \left( \frac{\rho_y}{\rho_z} \right)
\ee
and is a measure of similarity between these two random vectors.  This
divergence is nonnegative and is zero if and only if $\rho_z = \rho_y$
(\ie $Z_n  \egald Y_n$) \cite{CovTho91}.   This divergence has  also a
physical  signification, as  shown \eg  in  \cite{Qia01,FilHon05} Note
that KL divergence is not symmetric.  In the elliptical framework, its
expression can be simplified using the following lemma.

\begin{lemma}
  If  $Y_n$ and  $Z_n$  are elliptical  with  the same  characteristic
  matrix  $\Sigma_n=I_n$, then
\be
D_\kl(Y_n \, \| \, Z_n) = D_\kl(\| Y_n \| \, \| \, \| Z_n\|)
\label{KL_norm:eq}
\ee
\end{lemma}

\begin{proof}
By definition, 
\bens
D_\kl(Y_n \,  \| \,  Z_n) & =  & \int_{\mathbb{R}^{n}}  \rho_y(x) \log
\left( \frac{\rho_y (x)}{\rho_z (x)} \right) \d x
\\
& = &  \frac{2 \pi^{n/2}}{\Gamma(n/2)} \int_0^{+\infty} r^{n-1} d_y(r)
\log \left( \frac{d_y (r)}{d_z (r)} \right) \d r
\eens
with  $\rho_Y(x)  =  |  d_Y  \left(  \left(x^t  x\right)^{\frac{1}{2}}
\right)$ and  using \cite[4.642]{GraRyz80}.  But  using the expression
of $\rho_{\| Y_n \|}(r)$ as given by (\ref{Dn:eq}), we deduce
\bens
D_\kl(Y_n \,  \| \, Z_n)  = \int_0^{+\infty} \rho_{\| Y_n  \|}(r) \log
\left( \frac{\rho_{\| Y_n \|}(r)}{\rho_{\|  Z_n \|}(r)} \right) \d r =
D_\kl(\| Y_n \| \, \| \, \| Z_n\|)
\eens
\end{proof}
For  $\Sigma_n  \ne  I_n$,  $Y_n$   and  $Z_n$  must  be  replaced  by
$\Sigma^{-\frac{1}{2}}    Y_n$    and   $\Sigma^{-\frac{1}{2}}    Z_n$
respectively.

Applying this result  to the Student-t vector $X_n$  under study and a
zero-mean Gaussian vector $G_n$ with the same covariance matrix yields
the following results.

\begin{theorem}
  For $m  > 2$, the  KL divergences between  a Student-t vector  and a
  Gaussian vectors with the same covariance matrix are given by
\be\ba{lll}
D_\kl(X_n  \, \|  \,  G_n)  & =  &   \frac{n}{2} \log  \left(
  \frac{2 e}{m-2} \right) + \log \Gamma \left( \frac{n+m}{2} \right) -
\log \Gamma \left( \frac{m}{2} \right)\vspace{2.5mm}
\\
&  & +  \frac{m+n}{2} \left(  \psi \left(  \frac{m}{2} \right)  - \psi
  \left( \frac{n+m}{2} \right) \right)
\label{Dkl_XG_Student-t:eq}
\ea\ee
and
\be\ba{lll}
D_\kl(G_n  \, \| \,  X_n) &  = &  - \frac{n}{2}  \log \left(
  \frac{2 e}{m-2} \right) - \log \Gamma \left( \frac{n+m}{2} \right) +
\log \Gamma \left( \frac{m}{2} \right)\vspace{2.5mm}
\\
& & + \frac{m+n}{2} J(n)
\label{Dkl_GX_Student-t:eq}
\ea\ee
where
\be
J(n)  =  \frac{\displaystyle  \int_0^{+  \infty}  r^{n-1}  \log(1+r^2)
  \mbox{e}^{- \, \frac{m-2}{2}  \, r^2}}{2^{\frac{n}{2}-1} (m-2)^{- \,
    \frac{n}{2}} \Gamma \left( \frac{n}{2} \right)}
\label{Jn:eq}
\ee
Moreover\footnote{$f \sim g$ means that $f/g \to 1$ or $f = g + o(g)$}
when $n\to\infty$ and $m \to \infty$ with $n$,
\be
\hspace{-7mm} D_\kl(X_n \, \|  \, G_n) \sim \frac{n}{2} \log
\left( \frac{m}{m-2}  \right) + \frac{1}{2}  \log \left( \frac{m}{n+m}
\right) - \frac{n}{2 m} - \frac{n (n + 2 m)}{6 m^2 (n + m)}
\label{Dkl_XG_Student-t_asym:eq}
\ee
whereas, when $m = O(1)$ we have
\be
D_\kl(X_n  \, \| \,  G_n) =  \frac{n}{2} \left(  \log \left(
    \frac{2}{m-2} \right) + \psi  \left( \frac{m}{2} \right) \right) +
o(n)
\label{Dkl_XG_Student-t_asym_cte:eq}
\ee
and bounds for function $J(n)$ are given by
\be
J(n)  \:   \ge  \:  \psi\left(  \frac{n}{2}  \right)   -  \log  \left(
  \frac{m-2}{2}  \right) +  \log \left(  1 +  \frac{m-2}{2}  \, \left(
    \frac{\Gamma    \left(    \frac{n}{2}    \right)}{\Gamma    \left(
        \frac{n+1}{2} \right)} \right)^2 \right)
\label{Jn_lower:eq}
\ee
\be
J(n)  \:  \le \:  \log  \left( 1+  \frac{n}{m}  \right)  + \frac{2  \,
  n}{(n+m) \, (m-2)}
\label{Jn_upper:eq}
\ee
\end{theorem}

\begin{proof}
  Note first  that if  $m>2$ the covariance  matrix of $X_n$  with pdf
  (\ref{pdf_Student-t:eq}) is  $\frac{1}{m-2} I_n$ \cite{CosHer03}; if
  $m\le2$,  this  covariance  does  not  exist.   Then  starting  from
  (\ref{KL_norm:eq}),  using  (\ref{Dn_Student-t:eq})  for  $X_n$  and
  (\ref{fndn:eq})-(\ref{Dn:eq})    for    $G_n$    (with    covariance
  $\frac{1}{m-2}  I_n$), the  KL  divergence between  $X_n$ and  $G_n$
  writes
\bens\ba{l}
\hspace{-5mm} D_\kl(X_n \, \| \, G_n) = \frac{n}{2} \log 2 +
\log  \Gamma  \left(  \frac{n+m}{2}   \right)  -  \log  \Gamma  \left(
  \frac{m}{2} \right) - \frac{n}{2} \log(m-2)\vspace{2.5mm}
\\
\hspace{5mm} + \frac{m-2}{2} \displaystyle{\int_0^{+\infty} r^2 D_n(r)
  \,  \d r}  - \textstyle{\frac{m+n}{2}  \, \frac{2  \,  \Gamma \left(
      \frac{n+m}{2}  \right)}{\Gamma  \left(  \frac{n}{2}  \right)  \,
    \Gamma \left( \frac{m}{2} \right)}} \displaystyle{\int_0^{+\infty}
  \frac{r^{n-1}}{(1+r^2)^{\frac{n+m}{2}}} \, \log(1+r^2) \, \d r}
\ea\eens
The first integral  term is equal to $\frac{n}{m-2}$  while the second
one  is  evaluated  noticing  that  $\int h^k(r)  \log(h(r))  \d  r  =
\frac{\partial}{\partial    \lambda}    \left.    \int    h^\lambda(r)
\right|_{\lambda = k}$  and with \cite[8.380--3]{GraRyz80}, leading to
(\ref{Dkl_XG_Student-t:eq}).

The KL divergence between $G_n$ and $X_n$, (\ref{Dkl_GX_Student-t:eq})
is derived from the same technique.

Using  the  asymptotics  \cite[6.1.41  and  6.3.18]{AbrSte70}  of  the
log-gamma and of  the psi functions up to the  second order term, with
tedious algebra, (\ref{Dkl_XG_Student-t_asym:eq}) follows.

One can easily check that for any value $a$ and for $r \ge 0$ we have
\bens
\log(1+r^2) \le \log(1+a^2) + \frac{r^2 - a^2}{1+a^2}
\eens
Using  this   inequality  with  $a^2   =  n/m$  in  $J(n)$   leads  to
(\ref{Jn_upper:eq}).
Likewise, $ \log(1+r^2) = 2 \log r + \log(1+r^{-2})$ and for any value
$a$ and for $r \ge 0$ we have
\bens
\log(1+r^{-2}) \ge \log(1+a^{-2}) - \frac{2}{a (1+a^2)} (r-a)
\eens
Hence,    plugging   this   inequality    into   $J(n)$    and   using
\cite[4.352--1]{GraRyz80} leads to
\bens
J(n)  \ge  \textstyle \psi\left(  \frac{n}{2}  \right)  - \log  \left(
  \frac{m-2}{2} \right)  + \log ( 1  + a^{-2} )  - \frac{2}{a (1+a^2)}
\left(  \sqrt{\frac{2}{m-2}}   \,  \frac{\Gamma  \left(  \frac{n+1}{2}
    \right)}{\Gamma \left( \frac{n}{2} \right)} - a \right)
\eens
The best  bound is obtained  by maximizing the right-hand  side, which
amounts  to choose  $a =  \sqrt{\frac{2}{m-2}} \,  \frac{\Gamma \left(
    \frac{n+1}{2} \right)}{\Gamma \left( \frac{n}{2} \right)}$ leading
to (\ref{Jn_lower:eq}).
\end{proof}

Results from  the above theorem  show that $D_\kl$  does not
generally tend to zero with  $n$.  However, as previously noted, it is
more  significant  to  study  the  KL  divergence  rates  $\frac{1}{n}
D_\kl$.  Three situations occur:
\begin{itemize}
\item  If  $n  =  o(m)$, (\ref{Dkl_XG_Student-t_asym:eq})  shows  that
  $D_\kl(X_n \,  \| \, G_n)  $ tends to  0.  Thus $D_\kl(X_n \,  \| \,
  G_n) /  n$ tends also to  0 when $n$ increases.  This behavior shows
  that both in  the KL divergence sense and in  the KL divergence rate
  sense, a Gaussianization effect happens.
  Using the asymptotics  \cite[6.1.41 and 6.3.18]{AbrSte70}, we obtain
  that $D_\kl(G_n \, \| \, X_n) = \frac{n+m}{2} ( J(n) - \log(1+n/m) )
  +   \varepsilon(n)$   where    $\varepsilon(n)   =   o(1)$.    Using
  (\ref{Jn_upper:eq})  leads   to  $D_\kl(G_n   \,  \|  \,   X_n)  \le
  \frac{n}{m-2}  + \varepsilon(n)$.  Together  with the  positivity of
  the  KL  divergence, $D_\kl(G_n  \,  \| \,  X_n)$  tends  to 0  when
  $n\to\infty$: this confirms the  conclusion drawn from $D_\kl(X_n \,
  \| \, G_n)$.

\item     If     $m     \to     \infty$     in     $m     =     O(n)$,
  (\ref{Dkl_XG_Student-t_asym:eq})  tells us  that  the KL  divergence
  $D_\kl(X_n \, \|  \, G_n) $ has a finite non-zero  limit or can even
  diverge (\eg if  $m = o(n)$): again the  rate $\frac{1}{n} D_\kl(X_n
  \, \| \,  G_n)$ tends to 0 with $n \to  \infty$.  In this situation,
  there  is no  Gaussianization in  the  KL divergence  sense, but  in
  divergence rate the Gaussianization effect remains.
  From (\ref{Jn_lower:eq}),  the same technique  as in the  first case
  shows that the lower bound of $D_\kl(G_n  \, \| \, X_n) $ tends to a
  non-zero  limit  (and  clearly  diverges  if  $m=o(n)$).   From  the
  upperbound   (\ref{Jn_upper:eq}),  it   appears   that  $\frac{1}{n}
  D_\kl(G_n \,  \| \, X_n)$ tends  to 0 and the  conclusion drawn from
  $D_\kl(X_n \, \| \, G_n)$ holds.

\item If    $m   =   O(1)$    (\eg   $m$   constant),
\begin{itemize}
\item from (\ref{Dkl_XG_Student-t_asym_cte:eq})  one can conclude that
  no  Gaussianization appears,  neither in  KL divergence,  nor  in KL
  divergence rate.
\item from  the lower bound (\ref{Jn_lower:eq}), one  can check that
  $D_\kl(G_n  \,  \| \,  X_n)  $  diverges  with $n$.   However,  from
  (\ref{Jn_upper:eq}) it appears that  $\frac{1}{n} D_\kl(G_n \, \| \,
  X_n)$  tends  to 0.   This  behavior  seems  contradictory with  the
  previous one and tells that  in fact a Gaussianization effect exists
  in KL divergence rate.  This  contradiction is possible since the KL
  divergence  in  not  symmetric.   Note  also that  when  $m  \le  2$
  (\ref{Dkl_GX_Student-t:eq})-(\ref{Jn_lower:eq})-(\ref{Jn_upper:eq})
  can  still  be  considered   without  the  normalization  $m-2$  (no
  covariance for $X_n$): the conclusion holds in this case.
\end{itemize}
\end{itemize}

The  KL  divergence  rate   $\frac{1}{n}  D_\kl(X_n  \|  G_n)$  is  in
concordance with the observation  of the Gaussianization effect in the
distribution sense.   This generally holds  for $\frac{1}{n} D_\kl(G_n
\| X_n)$,  except notably when  $m$ does not  go to infinity:  in this
case,  although there  is no  Gaussianization in  distribution  the KL
divergence  rate goes to  0 with  $n$. Furthermore,  when $m$  goes to
infinity,  although  $X_n$ reaches  asymptotically  the  bound in  the
B.B.M.I., its distribution stays generally  at infinite -- or at least
at non zero distance -- ``distance'' (in the KL divergence sense) from
the Gaussian distribution.   Thus, much care must be  taken with these
conclusions, especially  about the real  meaning of the  KL divergence
rate.

\begin{theorem}
  For any $m > n-2$, the KL divergences between a Student-r vector and
  a Gaussian vector with the same covariance matrix is
\be\ba{lll}
D_\kl(X_n  \, \|  \,  G_n)  & =  &  \frac{n}{2} \log  \left(
  \frac{2 e}{m+2} \right) + \log \Gamma \left( \frac{m}{2} + 1 \right)
- \log \Gamma \left( \frac{m-n}{2} + 1 \right)\vspace{2.5mm}
\\
& &  + \frac{m-n}{2}  \left( \psi \left(  \frac{m-n}{2} + 1  \right) -
  \psi \left( \frac{m}{2} + 1 \right) \right)
\label{Dkl_XG_Student-r:eq}
\ea\ee

Moreover it verifies the asymptotic
\be
D_\kl(X_n \, \|  \, G_n) \: \sim \:  \frac{1}{2} \log \left(
  \frac{m+2}{m-n+2}   \right)    -   \frac{n   \,    (m-n)}{(m+2)   \,
  (m-n+2)}\nonumber
\label{Dkl_XG_Student-r_asym:eq}
\ee
\end{theorem}

\begin{proof}
  Note  first  that the  covariance  matrix  of  $X_n$ as  defined  by
  (\ref{pdf_Student-r:eq})  is  $\frac{1}{m+2}  I_n$  \cite{CosHer03}.
  Then,  (\ref{Dkl_XG_Student-r:eq}) is  obtained in  the same  way as
  (\ref{Dkl_XG_Student-t:eq}).

  Notice moreover  that here,  since $m \ge  n$, $m$ goes  to infinity
  with  $n$.   Using  ,  as  in the  Student-t  case,  the  asymptotic
  expansions \cite[6.1.41  and 6.3.18]{AbrSte70} of  the log-gamma and
  of   the   psi   functions    up   to   the   second   order   term,
  (\ref{Dkl_XG_Student-r_asym:eq}) follows.
\end{proof}

Thus, in the Student-r case, two behaviors arise:
\begin{itemize}
\item If $n  = o(m)$, from (\ref{Dkl_XG_Student-r_asym:eq}) $D_\kl(X_n
  \, \| \, G_n)$ tends to  zero $n$, hence the KL divergence rate goes
  also to 0:  both in the KL divergence and in  the KL divergence rate
  sense a Gaussianization effect appears.
\item If $n = O(m)$, the KL  divergence does not tends to zero and can
  even  diverge  with  $n$ if  $m-n  =  o(n)$  (e.g.  in  the  uniform
  case). But $D_\kl(X_n \, \| \, G_n)  / n$ still tends to 0 with $n$,
  exhibiting again an asymptotic Gaussianization behavior.
\end{itemize}

These observations  can be linked to  a famous result  by Diaconis and
Freedman   \cite{DiaFre87}  that   quantifies   the  total   variation
divergence\footnote{ We recall that  the total variation divergence is
  the $\L_1$-norm difference $D_\tv(Y  , Z) = \int_{\Rset^n} \| \rho_Y
  - \rho_Z \| $} between $X_n$ and $G_n$ as
\be
D_\tv(X_n,G_n) \le\frac{2 (n+3)}{m-n-1}  \hspace{5mm} \forall \: 1 \le
n \le m-2.
\ee
for integer $m$, where $X_n$ is built from the $n$ first components of
a $(m+2)$-dimensional  vector uniformly distributed in  the surface of
the $(m+2)$-dimensional sphere \cite{DiaFre87}.  Moreover, a necessary
and sufficient  condition for convergence to 0  of $D_\tv(X_n,G_n)$ is
$n=o(m)$.   These  results  were  extended  to the  KL  divergence  by
O. Johnson \cite{Joh03} as follows:
\be
D_\kl(X_n  \, \|  \,  G_n)  \le \log  \left(  \frac{m}{m-n} \right)  +
\frac{2}{\sqrt{m+2}/C-1}
\ee
A converse  result is also  provided in \cite{Joh03}, stating  that if
$D_\kl(X_n \, \| \, G_n)\rightarrow 0$ then $n=o(m)$.

%%%%%%%%%%%%%%%%%%%%%%%%%%%%%%%%%%%%%%%%%%%%%%%%%%%%%%%%%%%%

\section{Concluding remarks and discussion}
\label{conclusion:sec}

In  this paper, we  first extended  the entropic  uncertainty relation
found by Bialynicki-Birula \& Mycielski to R\'enyi entropies.  We have
checked that, for  a given dimension $n$, the  bound is again attained
in the  Gaussian case  and in this  case only. We  analytically proved
that the bound is also  asymptotically attained with the dimension $n$
in the conjugate multivariate  exponential power case, whatever $p \ge
1$.  We  numerically showed that the bound  is asymptotically attained
for the Cauchy  case for any value of $p>2$,  extending the results of
Abe, as  well as  in the general  Student-t case, including  these two
cases. This asymptotic analysis  was confirmed analytically.  The same
conclusion was drawn in the Student-r context, as far as our numerical
simulations are  valid.  These results  seem to violate the  fact that
the  bound  is  attained {\em  only}  for  the  Gaussian case.   If  a
Gaussianization effect was evoked in the first example, the second one
showed that  the effect of  dimension only may be  suspected, provided
some  favorable   conditions  exist.   To   get  a  feeling   to  this
interpretation,   let   us  come   back   to   the  Beckner   relation
(\ref{Beckner:eq}),  and assume  that  we deal  with  a wave  function
$\Psi_n$ such that
\be
(C_{p,q})^n \,  \frac{\|\Psi_n\|_p}{\|\widehat{\Psi}_n\|_q} = h(n,p) +
o(h(n,p))
\ee
for large  $n$ and  for some function  $h(n,p)$.  Again by  taking the
logarithm and  using the same approach as  in section \ref{Renyi:sec},
it is easy to show that
\bens
\hspace{-5mm}               \frac{H_{\frac{p}{2}}(X_n)               +
  H_{\frac{q}{2}}(\widetilde{X}_n)  }{n}  =  \log(2\pi)  +  \frac{\log
  p}{p-2} +  \frac{\log q}{q-2} +  \frac{2 p \log(h(n,p))}{n  (2-p)} +
o(1)
\eens
As a conclusion, it  is sufficient that $\frac{\log(h(n,p))}{n}$ tends
to  0 as  $n$ increases:  $h(n,p)$  does not  need to  converge to  1,
showing  that  no Gaussianization  effect  is  needed  here; in  fact,
$h(n,p)$  can even  diverge with  $n$.  This  conclusion holds  in the
Shannon  case, provided that  $\frac{h(n,p)}{p-2}$ (and  the remaining
term) has a limit as $p \to 2$.  In fact, in the conjugate exponential
power case (Student-t with  $m=n+2$), function $h(n,p)$ writes $h(n,p)
= h(p) = \frac{\log p}{2 \, p} - \frac{\log 2}{2 \, q} $, illustrating
this  conclusion.   As  perspective,  the  Student-r  case  should  be
analytically solved  to confirm  the numerical investigations.   To go
further,  it  seems interesting  to  determine  what  are the  minimal
conditions a pdf $f_n$ should verify so that it attains asymptotically
the bound  in the B.B.M.I.  As  far as we know,  this question remains
open  and we  still investigate  it. We  suspect however  that  in the
elliptical  context, no  major  additional constraints  are needed  to
reach the same conclusion.  Indeed, we feel that in this framework the
normalization term $1/n$  in the entropy rates may  be strong: one can
invoke Poincar\'e's observation  inducing a ``Gaussianization'' and we
suspect that the remaining contribution  of the sum of the entropy can
diverge, but at most in $o(n)$.

%%%%%%%%%%%%%%%%%%%%%%%%%%%%%%%%%%%%%%%%%%%%%%%%%%%%%%%%%%%%%%%%%%%%%

%%%%%%%%%%%%%%%%%%%%%%%%%%%%%%%%%%%%%%%%%%%%%%%%%%%%%%%%%%%%%%%%%%%%%

\bibliography{PhysicaA06}
\bibliographystyle{elsart-num}

\end{document}